\documentclass[preprint, 11pt]{elsarticle2}

\usepackage{mathrsfs}
\usepackage{amsfonts}
\usepackage{amsthm}
\usepackage{paralist}
\usepackage{graphics}
\usepackage{dsfont}
\usepackage{epstopdf}
\usepackage{float}
\usepackage{graphicx,subfigure}
\usepackage{amsmath}
\usepackage{bm}
\usepackage[colorlinks=true]{hyperref}
\usepackage{color}
\usepackage{epstopdf}
\hypersetup{urlcolor=blue, citecolor=red}
\usepackage{lineno,hyperref}

\newtheorem{theorem}{Theorem}[section]
\newtheorem{lemma}[theorem]{Lemma}

\theoremstyle{definition}
\newtheorem{definition}[theorem]{Definition}

\newtheorem{rem}{Remark}[section]
\theoremstyle{remark}

\newtheorem{cor}{Corollary}[section]

\numberwithin{equation}{section}

\def\eVO{e_{{V}}}
\def\eo{e_{\omega}}

\def\ub{\bm{u}}

\begin{document}

\begin{frontmatter}

\title{A highly efficient second-order long-time-dynamics-preserving scheme for geophysical fluid models}


\author[a]{Daozhi Han
		}
		\ead{daozhiha@buffalo.edu}

		\author [b]{Xiaoming Wang
           }
		\ead{wxm.math@outlook.com}

		\address[a]{Department of Mathematics, State University of New York at Buffalo, Buffalo, NY 14260}
		\address[b]{School of Mathematical Sciences, Eastern Institute of Technology, Ningbo, Zhejiang 315200, China}

\date{}

\begin{abstract}
	We develop and analyze a highly efficient, second-order time-marching scheme for infinite-dimensional nonlinear geophysical fluid models, designed to accurately approximate invariant measures—that is, the stationary statistical properties (or “climate”) of the underlying dynamical system. Beyond second-order accuracy in time, the scheme is particularly well suited for long-time simulations due to two key features: (i) it requires  solving only a fixed symmetric positive-definite linear system with constant coefficients  twice
    at each step, and (ii) it guarantees {  unconditional} long-time stability, producing uniformly bounded solutions in time for any bounded external forcing, regardless of initial data. For prototypical {  damped-forced} models such as the {  damped-forced continuously stratified} quasi-geostrophic equation, the method preserves dissipativity, ensuring that numerical solutions remain bounded in a function space compactly embedded in the phase space as time tends to infinity. Leveraging this property, we rigorously prove convergence of both global attractors and invariant measures of the discrete system to those of the continuous model in the vanishing time-step limit.
	
	A central innovation of the method is a mean-reverting scalar auxiliary variable (mr-SAV) formulation, which preserves the dissipative structure of externally forced systems within an appropriate phase space. For the infinite-dimensional models considered, we additionally employ fractional-order function spaces to ensure the well-posedness of the scheme within the same space and to establish { asymptotic} compactness of numerical solutions in topologies compatible with the phase space.

    To the best of our knowledge, this is the first second-order scheme that simultaneously (a) requires two solves of a fixed symmetric positive-definite system 
    at each step, (b) inherits dissipativity of the infinite-dimensional nonlinear model without step-size restrictions, and (c) provably captures the long-time statistical behavior of the underlying system.
{ Numerical results on the 2D Navier–Stokes equations demonstrate the accuracy, unconditional long-time stability, and robustness of the scheme. Moreover, long-time simulations indicate that several statistical indicators of the maximum vorticity—such as the post-transient mean, variance, and burst exceedance probabilities—exhibit a degree of consistency under successive temporal refinements. This provides numerical evidence that the scheme captures key qualitative features of the long-time dynamics, including intermittent burst behavior, in a manner that is consistent across temporal refinements.}
    
\end{abstract}

\end{frontmatter}

	\section{Introduction} \label{intro-sec}

Many complex systems display chaotic behavior with sensitive dependence on initial data, making it impossible to achieve precise long-term prediction. It is well known that the statistical properties of such systems are more physically relevant \cite{LaMa1994, Frisch1995, FMRT2001, SiSt2001, MaWa2006a, MoYa2007}. In particular, if a complex system attains statistical equilibrium and the long-time statistics are of interest, it is the global attractor and the invariant probability measure carried by it that are important. A classical example is the climate system, which is a forced dissipative chaotic dynamical system. Mathematically, climate is defined as a statistical ensemble of states assumed by the climate system over a sufficiently long time interval---30 years according to the World Meteorological Organization. The physically relevant characteristic of climate is thus the essential invariant measure concentrated on the attractor \cite{DyFi1997, Majda2012a}. Another example is the spontaneous occurrence of extreme events, defined by abrupt changes in the state of a dynamical system of many degrees of freedom, such as oceanic rogue waves, extreme weather patterns, earthquakes, shocks in power grids, and stock market crashes, among many others \cite{FaSa2019}. The interests again are on the complex high-dimensional attractor and the invariant measure supported on the attractor. Since direct physical experiments are generally not feasible, numerical experiments become essential for these studies. It is therefore highly desirable to develop efficient numerical methods capable of capturing the long-term dynamics characterized by global attractors and invariant measures of the underlying dissipative dynamical systems.

Long-time numerical integration of dynamical systems always poses a significant challenge since small errors could accumulate over long time intervals unless the dynamics are asymptotically stable (no chaos) \cite{HeRa1986}. For chaotic dynamical systems with positive Lyapunov exponents, individual trajectories cannot be accurately approximated over long time scales \cite{StHu1996}. {Nevertheless, the numerical approximation of long-time statistics, long understood as a central issue \cite{Stuart1994, SiSt2001}, is possible, as demonstrated in one of the authors’ works, at least for first-order time-marching methods \cite{Wang2010a, Wang2016}}. Aspects of long-time dynamics of time-marching schemes are systematically treated by Stuart and Humphries in the book \cite{StHu1996}. Elliott and Stuart in \cite{ElSt1993} examined both explicit and implicit schemes from the perspective of preserving the steady-state solution to a class of semi-linear parabolic equations of gradient type (see also \cite{HJS1994, JST1998}). Hale et al. in \cite{HLR1988} identified sufficient conditions ensuring convergence of attractors in the upper semi-continuity fashion for a one-step scheme, which was then strengthened by Hill and S\"{u}li in \cite{HiSu1996}. Hill and S\"{u}li further established corresponding results for linear multistep methods in \cite{HiSu1995} and for the 2D Navier--Stokes equations (NSE) with a semi-explicit treatment of the nonlinear advection term in \cite{HiSu2000}. For dissipative dynamical systems generated by continuous semi-groups, Wang in \cite{Wang2010a, Wang2016} established Lax-type criteria for the convergence of invariant measures (stationary statistical solutions) of a one-step scheme. Application of the general framework to the 2D Rayleigh--B\'enard convection was carried out by Tone and Wang in \cite{ToWa2011}, to double diffusive convection by Tone et al. in \cite{TWW2015}, and to the 2D Navier--Stokes equations in the vorticity--stream function formulation by Gottlieb et al. in \cite{GTWWW2012}. Separately, in \cite{Wang2012}, Wang established convergence of long-time statistical properties of a BDF2 time-marching scheme at vanishing time step for the 2D Navier--Stokes equations with periodic boundary conditions. All these works focus on time discretizations with explicit treatment of the nonlinear advection/convection terms under a time-step constraint proportional to the inverse of the Reynolds number. 

In recent years, a Lagrange-multiplier-type approach, termed the Scalar Auxiliary Variable (SAV) method, was first developed by Shen et al. for solving PDEs with a gradient flow structure \cite{SXY2018, SXY2019}, see also \cite{GuTi2013, YZW2017, YaJu2017, YaHa2017} for early variants of this method. High-order SAV methods are constructed in \cite{ALL2019, GZW2020, HuSh2021, HuSh2022, WHS2022}. More recently, the SAV approach has been extended to treat the Navier--Stokes equations (NSE) and related fluid models; see, for instance, \cite{LYD2019, LiSh2020, LSL2021, Yang2021, JiYa2021, CHJ2023, DHJL2025} among many others. {The extension to NSE-type models is termed ZEC or extended-ZEC by Yang \cite{Yang2021}. See also \cite{LSL2021}.} In particular, the advective term of the NSE is treated explicitly in an SAV { (ZEC)} method, and a prefactor (a scalar auxiliary variable) is discretized implicitly. As a result, the SAV methods lead to the same Stokes solver across the time iterations, yielding a highly efficient algorithm for long-time computation. In the absence of external forcing, it was shown that the SAV { (ZEC)} methods for the NSE are unconditionally stable in the sense that a (modified) energy is non-increasing \cite{LSL2021} and the references therein. In the case of a linear advection--diffusion equation with constant coefficients, the stability of the SAV method can be attributed to the reduction of the effective P\'eclet number. For fixed final time and with external forcing, the SAV methods designed in \cite{LiSh2020, Yang2021} are also unconditionally stable. To the best of our knowledge, there are no existing SAV { (ZEC)} methods for the NSE that can preserve uniform-in-time estimates of the velocity in the energy norm and long-time estimates in the $H^1$ norm in the 2D case in the presence of nontrivial external forcing. Recently, in \cite{CHW2024}, we developed a mean-reverting-SAV-BDF2 scheme for autonomous ODE systems with skew-symmetric nonlinearity reminiscent of the Navier--Stokes equations. It was shown that the discrete global attractor and invariant measure both converge to their continuous counterparts at vanishing time step. 

In this article, we extend the \textbf{mean-reverting SAV-BDF2} method to infinite-dimensional dynamical systems and investigate the long-time dynamics of the scheme for solving geophysical fluid models, including the 2D damped-driven barotropic quasi-geostrophic (QG) model, the 2D NSE in vorticity--stream function formulation, and the 3D damped-driven continuously stratified QG equations. {   A major obstacle in the infinite-dimensional setting is the lack of compactness of bounded subsets of the phase space. Consequently, proving the existence and convergence of invariant measures requires large-time stability bounds 
for the numerical solution in function spaces that are compactly embedded into the phase space. Deriving these higher-regularity estimates is highly nontrivial due to the presence of unbounded differential operators and 
nonlinear advection terms.} 

For external forcing bounded uniformly in time, a uniform-in-time estimate of the numerical solution in the $L^2$ norm is established for any time-step size for the general scheme. {Unlike previous SAV methods, our equation for the auxiliary scalar variable is mean-reverting with a constant rate. The constant-rate mean-reverting nature of the auxiliary variable equation dictates that the auxiliary variable approaches the desired value of one exponentially fast regardless of the initial data. This yields asymptotic equivalence between the expanded system and the original model at large time for all initial values.} For external forcing independent of time, this freedom of initial condition enables us to write the scheme as a discrete autonomous dynamical system on the product space $(H^\alpha \times \mathbb{R})^2$, { $\alpha\in [\frac{1}{4}, 1)$, building on the dynamical system approach embedded in the Lax-type criteria} \cite{Wang2010a, Wang2016}.  {  Crucially, we  establish a uniform 
$H^1$ estimate of the numerical solution at large time, which combined with the compact embedding $H^1 \hookrightarrow H^\alpha$ implies that the discrete dynamical system admits a global attractor uniformly bounded in the $H^1$ norm with respect to the time-step sizes}.   By focusing on the global attractor, and leveraging asymptotic consistency and finite-time trajectory convergence {  uniformly with respect to} $H^1$ initial data, we establish convergence of the global attractors and invariant measures (stationary statistical solutions) of the numerical scheme to those of the continuous models at vanishing time step. 

{  The mean-reverting term can be understood as a damping term, and similar damping was utilized in some of the extended ZEC work by Yang and collaborators \cite{ZHY2025}. However, it is not clear whether the solution-dependent nonlinear damping utilized in the extended ZEC schemes  \cite{ZHY2025} 
can guarantee uniform-in-time bound on the solution in the presence of a non-trivial forcing.} As far as we know, the novel mean-reverting-SAV-BDF2 method is the only algorithm that treats the nonlinear advection term in the above-mentioned geophysical fluid models (essentially) explicitly while preserving uniform-in-time estimates without any time-step restriction in the presence of external forcing. We refer to \cite{SiAr1994, HeRa1990, ToWi2006, Tone2007, ChWa2016, HOR2017, ReTo2023, Shiue2025} for the design of long-time $H^1$-stable schemes for solving the 2D NSE.

{ 
\subsection{The main results}

The major theoretical results of this article include
\begin{enumerate}
    \item unconditional time-uniform estimate of the numerical solution in the $L^2$ norm, Theorem \ref{theo1};
    \item large-time uniform $H^1$ bound of the scheme for the 2D NSE and the 3D continuously stratified QG model, Theorem \ref{main-co};
    \item existence of global attractors and invariant measures generated by the scheme, Corollary \ref{Uni-Cor};
    \item convergence of global attractors and invariant measures, Theorem \ref{att-conv} and Theorem \ref{im-conv}. 
\end{enumerate}
}

{ The proof has four main ingredients. First, the BDF2 G-stability and the skew-symmetry of the nonlinear term yield an unconditional absorbing estimate in the base Hilbert space. Second, for the continuously stratified QG model, a large-time $H^1$ estimate provides asymptotic compactness in $H^\alpha, \alpha \in [\frac{1}{4}, 1)$. Third, asymptotic consistency shows that the two components of the BDF2 phase vector become close and that the auxiliary variable approaches one. Finally, these estimates are combined with finite-time trajectory convergence to prove upper semi-continuity of attractors and convergence of invariant measures.}

{The results in this article provide significant improvement over previous results on long-time dynamics of systems like the 2D NSE \cite{HiSu2000, Wang2012} in several aspects.} First, the BDF2 scheme in \cite{HiSu2000} employed a semi-implicit discretization of the advection term, resulting in a coefficient matrix varying in time. In contrast, the mean-reverting SAV-BDF2 scheme proposed in this article treats the advection term essentially fully explicitly, leading to the same coefficient matrix across all time levels and is therefore highly efficient for long-time computation. {Secondly, only the global attractor was treated in \cite{HiSu2000} without referring to the invariant measures.} Thirdly, there are time-step size constraints in \cite{Wang2012} due to explicit discretization of the nonlinear term, and in \cite{HiSu2000} to obtain the $H^1$ uniform boundedness of the global attractor, whereas the uniform $H^1$ estimate in our work is obtained unconditionally. Fourthly, \cite{HiSu2000} utilized the monoid approach from their earlier work \cite{HiSu1995} to analyze convergence, which requires $H^1$ initial data. Hence, the existence of a global attractor is achieved for the fully discrete scheme, in the absence of uniform $H^2$ estimates. {Fifthly, \cite{Wang2012, ChWa2016} dealt exclusively with periodic boundary conditions, while the current work applies to a more general setting.} On the other hand, we adopt a dynamical system formalism on the product space $(H^\alpha \times \mathbb{R})^2$, with a dimension dependent $\alpha$, allowing us to establish the existence of a global attractor for the semi-discrete scheme. Finally, we establish the convergence of both the attractors and the stationary statistical solutions (invariant measures), while \cite{HiSu2000} only dealt with attractor convergence and \cite{Wang2012} only established convergence of invariant measures. Finite-time numerical approximations of statistical solutions in various settings and their convergence analyses are reported in several recent works \cite{MiSc2012, LMS2016b, FLM2018, FLMW2020, LMP2021, Par2023}.

The rest of the article is organized as follows. { In Sec. \ref{sec(2)} we introduce the general mean-reverting SAV-BDF2 scheme  applicable to the continuously stratified quasi-geostrophic equation. 
We then  establish the unconditional uniform-in-time bound {on the solution to} the scheme in the phase space in Sec. \ref{UIT}. In Sec. \ref{sta-H1} we further show  that for large time the numerical solution is also uniformly bounded in the $H^1$ norm. In Sec. \ref{Asym-sim}, we show that the solution to the numerical scheme is consistent with the continuous model in the sense that the scalar auxiliary variable is close to one and the difference of the numerical solution at two neighboring steps is small for large time.  
 The convergence of global attractors and invariant measures is then established in Sec. \ref{long-time}. 
 Numerical results related to accuracy, stability and robustness of the proposed methods are reported in Sec. \ref{s-nu}.  Some concluding remarks are provided in Sec. \ref{conc}. The Appendix section provides the finite-time convergence of the numerical scheme together with the proof of an auxiliary lemma. The finite time convergence is uniform for initial data from the global attractors of the numerical scheme.}




\section{The physical models and the scheme}\label{sec(2)}

{ 
In this section we will introduce the mr-SAV-BDF2 scheme that is applicable to a few damped-driven geophysical models.
We will first introduce the applicable physical models to fix the ideas. We then introduce the general scheme applicable to these models together with its efficient implementation that involves solving a fixed linear positive definite constant coefficient problem twice at each time step. 
}
\subsection{The physical models}
\subsubsection{Abstract formulation}
With applications in geophysical fluid models in mind, we consider the following forced dissipative infinite dimensional dynamical system in the form  
{ 
\begin{align}
	\frac{d\ub(t)}{dt}+{ \nu}A\ub(t) +N\big(\ub(t)\big)=\mathbf{F}(t), \quad \ub(t) \in \mathbf{H}, \label{DSe}
\end{align}
}
where   $A$ is a positive definite, self-adjoint  linear operator with compact resolvent on the Hilbert space $\mathbf{H}$, $\nu>0$ is a constant. Denote the inner product and the norm on $\mathbf{H}$ by $(\cdot, \cdot)$ and $||\cdot||$, respectively. Let $\mathbf{V}=\mathcal{D}(A^\frac{1}{2})$, and assume that
there is a positive constant $c_0$ (the first eigenvalue of $A$) such that 
\begin{align}\label{coer-es}
	||A^{\frac{1}{2}}\ub||^2\geq c_0||\ub||^2, \quad \forall \ub \in \mathbf{V}. 
\end{align}    
One has $\mathbf{V} \subset \mathbf{H}  \subset \mathbf{V}^\prime$ with continuous embeddings. The duality between $\mathbf{V}^\prime$ and $\mathbf{V}$ is identified as
$\langle \mathbf{f}, \mathbf{v} \rangle=(A^{-\frac{1}{2}}\mathbf{f}, A^{\frac{1}{2}}\mathbf{v}),  \mathbf{f} \in \mathbf{V}^\prime, \mathbf{v}\in \mathbf{V}$.
The nonlinear operator $N: \mathbf{V}  \rightarrow \mathbf{V}^\prime$ is assumed to be skew-symmetric, i.e. $\langle N(\ub), \ub\rangle=0$ for all $\ub\in \mathbf{V}$.  Throughout it is assumed that the body force $\mathbf{F}: (0, \infty)\rightarrow \mathbf{H}$ is uniformly bounded in time. It follows immediately from the dissipativity of the system that $\sup_{t>0}||\ub(t)|| \leq C$. { We note that additional assumptions are needed for each physical model and for asymptotic compactness arguments which we do not outline here; see \cite{Temam1997}.}

\subsubsection{The damped-driven continuously stratified quasi-geostrophic model}
Consider the  damped-driven continuously stratified QG equations \cite{MaWa2006a} in a periodic channel
$
\Omega=[0, 2\pi] \times [0, 2\pi]\times [0, 1]:
$
\begin{align}
	&\frac{\partial \omega}{\partial t}-D(\omega)  +\nabla^{\perp}_H \psi \cdot \nabla_H \omega +\beta \frac{\partial \psi}{\partial x}=f, \quad \omega =-\big(\Delta_H \psi+F^2 \frac{\partial^2 \psi}{\partial z^2}\big), \label{SVe}
\end{align}
where $\omega$ is the vorticity; $\psi$ is the stream function; { $\Delta_H=\partial_x^2+\partial_y^2$, $\nabla^{\perp}_H=[\partial_y, -\partial_x]^T$, $\nabla_H=[\partial_x, \partial_y]^T$  are the horizontal Laplacian,  perpendicular gradient and gradient operators, respectively}; $D(\omega):=\nu_H \Delta_H \omega+\nu_v \frac{\partial^2 \omega}{\partial z^2}$ is an anisotropic dissipation operator with the viscosity coefficients $\nu_H$ and $\nu_v$; $F>0$ is the Froude number which is the ratio of vertical to buoyancy restoring forces; $\beta$ is a constant representing the beta-plane effect from the Coriolis force. 
In addition to the horizontal periodic boundary conditions, in the $z$-direction one imposes either periodic BC or {homogeneous Dirichlet/Neumann BC
\begin{align}\label{BC}
	\omega=\psi=0, \text{ or } \frac{\partial \omega}{\partial z}=\frac{\partial\psi}{\partial z}=0, \quad z=0, 1.
\end{align}
For definiteness, we focus on the Dirichlet case, and introduce
\begin{align*}
	H_{p0}^1(\Omega)=\{\omega\in H^1(\Omega), \omega \text{ is periodic in $x$ and $y$ direction}, \omega|_{z=0, 1}=0\}.
\end{align*}
The functional setting is then $H=L^2(\Omega)$ and $V=H_{p0}^1(\Omega)$.
The model generates a dissipative dynamical system on $H$ \cite{Wang1992}.} 

\subsubsection{The barotropic quasi-geostrophic equations and the 2D Navier-Stokes equations}
The barotropic QG model (BQG) is a simplified two-dimensional system describing large-scale flows for the atmosphere or mesoscale flow in the oceans. Consider the damped-driven QG system in either a periodic domain or a two-dimensional periodic channel $\Omega=[0, 2\pi]\times [0, 1]$:
\begin{align}
  &\frac{\partial \omega}{\partial t}-\nu \Delta \omega  + \nabla^{\perp} \psi \cdot \nabla \omega+\beta \frac{\partial \psi}{\partial x}=F, \quad \omega=-\Delta \psi, \label{QG1} \\
  & \omega, \psi  \text{ periodic in } x, \quad \omega=\psi=0, \text{ on } z=0, 1. \label{QG2}
\end{align}
Here  $\omega$ is the vorticity, $\psi$ is the stream function,  and {$\nu$ is the kinematic viscosity}.   
The function spaces are $H=L^2(\Omega)$ and $V=H_{p0}^1(\Omega)$ which is a subspace of $H^1(\Omega)$ whose members are $2\pi$ periodic in the $x$ direction and vanish at $z=0,1$. 
The  model generates a dynamical system  \cite{Temam1997}.

When $\beta=0$ in \eqref{QG1}, the BQG system gives rise to the two dimensional Navier-Stokes equations (NSE) in vorticity-stream function formulation.

{ 
  The estimates below are presented for the continuously stratified QG model. It is clear that the BQG model  and the 2D NSE in vorticity–stream function formulation  are obtained as special cases of the same vorticity–stream function structure. In these cases the stream function elliptic relation, skew-symmetry of the transport term, Sobolev estimates, Poincaré-type inequalities, and compact embeddings used in Sections 3--6 remain valid, with only notational simplifications. Therefore the large-time $H^1$ bound, attractor convergence, and invariant-measure convergence arguments carry over to these periodic cases.
 
 }

\subsection{The general scheme and its efficient implementation}

\subsubsection{The extended model with mean-reverting SAV}
In order to formulate our scheme, we introduce a scalar auxiliary variable (SAV) $q(t)$ together with the following extended model (cf. \cite{CHW2024} for the case in finite-dimension)
\begin{eqnarray}
	\frac{d\ub (t)}{dt}+{ \nu}A\ub(t) +q\,N(\ub(t))&=&\mathbf{F}(t),\label{SAVe0}\\
	\frac{dq}{dt} -\langle N(\ub),\ub\rangle&=&-\gamma q+\gamma, \label{SAVe}
\end{eqnarray}
with a user-specified  positive constant parameter $\gamma$.
It follows { from the skew-symmetry of the nonlinear term} that 
\[q(t)=q_0 e^{-\gamma t}+(1-e^{-\gamma t})\rightarrow 1, \quad { t\rightarrow +\infty}.\] 
{  In particular, $q(t)\equiv1$ if $q_0=1$.}
\subsubsection{The mean-reverting SAV-BDF2 scheme}
Let $k>0$ be the time-step size, $t^n=n k$ for an integer $n$, and denote the numerical approximation of $\phi$ at $t^n$ by $\phi^n$. The second order backward differentiation (BDF2) approximation of the derivative is then defined as { 
\begin{align}
	\delta_t \phi^{n+1}:=\frac{3\phi^{n+1}-4\phi^n+\phi^{n-1}}{2k}. \label{BDF2-de}	
\end{align}	
Furthermore, {  the 2nd order extrapolation approximation} is denoted by $$\overline{\phi}^{n+1}:= 2\phi^n-\phi^{n-1}.$$}
The {\bf mean-reverting-SAV-BDF2 time-marching scheme (mr-SAV-BDF2)} based on the expanded system \eqref{SAVe0} and \eqref{SAVe} is as follows
\begin{align}
	&\delta_t \ub^{n+1} +{\color{blue}\nu}A  \ub^{n+1}+q^{n+1} N\big(\overline{\ub}^{n+1}\big) =\mathbf{F}^{n+1}, \label{DDSe0}\\
	&\delta_t q^{n+1}+\gamma q^{n+1}-\langle N\big(\overline{\ub}^{n+1}\big),  \ub^{n+1}\rangle=\gamma, \label{DSAVe0}
\end{align}	
{ where $\mathbf{F}^{n+1}=\mathbf{F}(t^{n+1})$.}
\subsubsection{Efficient implementation of the scheme}
The efficiency of the mr-SAV-BDF2 scheme is explained as follows. In light of the linearity of the scheme, one performs a linear superposition by setting
\begin{align}\label{Decomp}
	\ub^{n+1}=\ub_1 + q^{n+1} \ub_2,
\end{align}
where $\ub_i, i=1,2$ satisfy the equations
\begin{align}\label{Dec-sol}
	\frac{3 \ub_i}{2k}+ {\color{blue}\nu}A \ub_i=\mathbf{f}_i, \quad \text{ with } \mathbf{f}_1:=\mathbf{F}^{n+1}+	\frac{4 \ub^n-\ub^{n-1}}{2k}, \quad \mathbf{f}_2:=-N\big(\overline{\ub}^{n+1}\big).
\end{align}
{ 
Once \(u_1\) and \(u_2\) have been computed from \eqref{Dec-sol}, the scalar
auxiliary variable \(q^{n+1}\) can be obtained explicitly from the linear equation \eqref{DSAVe0}. Indeed,  we find
\[
        \left(
        \frac{3}{2k}
        +\gamma
        -
        \left\langle N(\bar u^{n+1}),u_2\right\rangle
        \right)q^{n+1}
        =
        \gamma
        +
        \frac{4q^n-q^{n-1}}{2k}
        +
        \left\langle N(\bar u^{n+1}),u_1\right\rangle .
\]
Moreover, since \(u_2\) satisfies
\[
        \frac{3}{2k}u_2+{ \nu}Au_2=-N(\bar u^{n+1}),
\]
hence
\[
        \frac{3}{2k}\|u_2\|^2+{ \nu}\|A^{1/2}u_2\|^2
        =
        -\left\langle N(\bar u^{n+1}),u_2\right\rangle .
\]
Therefore
\[
        q^{n+1}
        =
        \frac{
        \displaystyle
        \gamma
        +
        \frac{4q^n-q^{n-1}}{2k}
        +
        \left\langle N(\bar u^{n+1}),u_1\right\rangle
        }{
        \displaystyle
        \frac{3}{2k}
        +\gamma
        +
        \frac{3}{2k}\|u_2\|^2
        +
        { \nu}\|A^{1/2}u_2\|^2
        } .
\]
}
 The scheme is highly efficient and suitable for long-time computation because, after spatial discretization, it requires only linear solvers with a symmetric positive-definite coefficient matrix that does not change in time.


\section{Uniform-in-time estimate in $\mathbf{H}$}\label{UIT}
In this section, we establish the unconditional uniform-in-time estimate in $\bf{H}$ of the general mr-SAV-BDF2 scheme \eqref{DDSe0}--\eqref{DSAVe0}.
{ Following Dahlquist's G-stability framework \cite{Dahlquist1978, Hill1997}, one defines the symmetric positive definite $G$ matrix for the BDF2 scheme
\begin{align}
	G:=\frac{1}{4}\begin{bmatrix}
		5 & -2\\
		-2 & 1
	\end{bmatrix},
\end{align}
and the associated  $G$-norm on the product Hilbert space $\mathcal{H}$: $$||{V}||_G^2:=\big({V}, G{V}\big)=\frac{5}{4}||v_1||_{H}^2-(v_1, v_2)_H+\frac{1}{4}||v_2||_H^2$$ where ${V}=[v_{1}, v_2]^T \in \mathcal{H}=H\times H$. The following BDF2-G identity, cf. \cite{ReTo2023} and references therein, is useful
\begin{align}\label{G-identity}
    k\big(\delta_t u^{n+1}, u^{n+1}\big)_H=||V^{n+1}||_G^2-||V^{n}||_G^2+\frac{1}{4} ||u^{n+1}-2u^n+u^{n-1}||_H^2,
\end{align}
where $V^{n+1}=[u^{n+1}, u^n]^T$.
The $G$-norms and the standard norms are equivalent in the sense that
\begin{align}
	&c_l ||{V}||_G^2 \leq || {V}||_{\mathcal{H}}^2 \leq c_u ||{V}||_G^2, \label{equi1}
\end{align}
where $c_l, c_u$ are the inverses of  the largest and  smallest eigenvalue of $G$, respectively. In the sequel, the subscript $H$ will be omitted if it is self-evident from the context. The $G$-norm is defined similarly for $\mathcal{H}=\mathbb{R}\times \mathbb{R}$, denoted by $|\cdot|_G$.}

{ 
The following lemma is crucial to the dissipative estimates of the BDF2 scheme. The proof  is deferred to the Appendix I.

\begin{lemma}\label{Lem-key}
Let $k>0$, $g\ge 0$, $l_k:=\min\{k,1\}$, and let $M\ge 1$ be an integer.
Let
\[
\{Y_n\}_{n\ge M},\qquad \{r_n\}_{n\ge M-1}
\]
be nonnegative sequences. Assume that there exist  constants $c_d>0$,  $\mu>0$ and $0\le \varepsilon\le 1/12$ such that
for all $n\ge M$,
\begin{align}
r_n+r_{n-1}&\ge c_dY_n,
\label{Bas_as} \\
Y_{n+1}-Y_n+\mu k r_{n+1}
&\le
\varepsilon\mu k(r_n+r_{n-1})+kg.
\label{Bas_ine}
\end{align}

Then, for any $0<b\le \min\left\{1,\frac{\mu c_d}{24}\right\}$ and for all $n\ge M$,
\begin{equation}
Y_{n+1}
\le \max\left\{1,\frac{48}{\mu c_d}\right\} \Big[
\frac{
Y_M+\mu \left(\frac14+\varepsilon\right)r_M
+\varepsilon\mu  r_{M-1}
}{
(1+bl_k)^{\,n+1-M}
}
+
\frac{g}{b}\Big].
\label{lem_fr}
\end{equation}

\end{lemma}

}

Recall that $V^{n}:=[\ub^{n}, \ub^{n-1}]^T$, $Q^{n}:=[q^{n}, q^{n-1}]^T$.  The unconditional uniform-in-time estimate in $\bf{H}$ is established in the following Theorem.
\begin{theorem}\label{theo1}
{Let \(\big(\mathbf{u}^i, q^i\big) \in \mathbf{V}\times \mathbb{R},  i=0, 1\)},  \(k>0\), \(l_k:=\min\{k,1\}\), and assume that
\(F\in L^\infty(0,\infty;\mathbf{H})\).
Then there exist a positive constant $C_g$ and an integer $N_0(k)$ such that the solution of the mr-SAV-BDF2 scheme satisfies
\begin{align}\label{zero-es}
   \|V^{n}\|_G^2+|Q^{n}|_G^2
        \le R_0^2, \quad \forall n\ge N_0,
\end{align}  
 where 
        \[R_0^2:=C_g
        \left(
        \frac{1}{c_0 \nu}\|F\|_{L^\infty(0,\infty;H)}^2+\gamma
        \right).  \]
\end{theorem}

\begin{proof}

Testing Eq. \eqref{DDSe0} with $\ub^{n+1} k$, multiplying Eq. \eqref{DSAVe0} by $q^{n+1} k$, and adding the two resulting identities, the
nonlinear terms cancel. Using the BDF2 \(G\)-identity \eqref{G-identity}, we obtain

\begin{align}\label{bas-l2}
        &\|V^{n+1}\|_G^2-\|V^n\|_G^2
        +
        |Q^{n+1}|_G^2-|Q^n|_G^2  +                                    
        \frac14\|\mathbf u^{n+1}-2\mathbf u^n+\mathbf u^{n-1}\|^2 \nonumber \\
        &
        \quad +
        \frac14|q^{n+1}-2q^n+q^{n-1}|^2                               
        +
        k{\color{blue}\nu}\|A^{1/2}\mathbf u^{n+1}\|^2
        +
        \gamma k|q^{n+1}|^2              \nonumber                               \\
        &=
        k(F^{n+1},\mathbf u^{n+1})
        +
        \gamma k q^{n+1}.
\end{align}
By the Cauchy-Schwarz inequality and Poincar\'e-type inequality,
\[
        k(F^{n+1},\mathbf u^{n+1})
        \le
        \frac{c_0 \nu k}{2}\|\mathbf u^{n+1}\|^2
        +
        \frac{k}{2c_0 {\color{blue}\nu}}\|F^{n+1}\|^2,\quad 
        \gamma k q^{n+1}
        \le
        \frac{\gamma}{2}k|q^{n+1}|^2
        +
        \frac{\gamma}{2}k .
\]
By the assumption \eqref{coer-es},  we get
\[
        \|V^{n+1}\|_G^2+|Q^{n+1}|_G^2
        +
        \frac{c_0 \nu}{2}k\|\mathbf u^{n+1}\|^2
        +
        \frac{\gamma}{2}k|q^{n+1}|^2
        \le
        \|V^n\|_G^2+|Q^n|_G^2
        +
        \frac{k}{2c_0{\color{blue}\nu}}\|F^{n+1}\|^2
        +
        \frac{\gamma}{2}k .
\]

Define
\[
Y_n:=\|V^{n}\|_G^2+|Q^{n}|_G^2, \quad r_n:=\|\mathbf u^{n}\|^2+|q^{n}|^2, 
\]
and
\[
g:=\frac{1}{2}
        \left(
        \frac{1}{c_0 \nu}\|F\|_{L^\infty(0,\infty;H)}^2+\gamma
        \right), \quad \mu:=\frac12\min\{c_0 \nu,\gamma\}.
\]
One obtains
\begin{align}\label{Th-es-1}
         Y_{n+1}+ \mu k r_{n+1} \leq Y_n+k g, \qquad n\geq 1,  
\end{align}
conforming to \eqref{Bas_ine} with $\epsilon=0$ and $M=1$. The assumption \eqref{Bas_as} holds  with $c_d:=c_l$ by  the norm equivalence \eqref{equi1}.

An application of Lemma \ref{Lem-key} then gives
\begin{align}
    \label{Fi1}
        Y_{n+1}
&\le \frac{C_0
}{
(1+b_0l_k)^{\,n}
}Y_1+ C_g g, \qquad n\ge 1, 
\end{align}
with
\[b_0:=\min\left\{1,\frac{\mu c_d}{24}\right\}, \quad C_0:=\max\left\{1,\frac{48}{\mu c_d}\right\} \big[1+c_u \mu/4\big], \qquad C_g:=\frac{1}{b_0}\max\left\{1,\frac{48}{\mu c_d}\right\}.\]
Noting that \(0<b_0l_k<1\) and \(2^x \leq 1+x, \forall x\in[0,1]\), one has
\[
\frac{C_0
}{
(1+b_0l_k)^{\,n-1}
}Y_1 \leq \frac{C_0
}{
2^{\,(n-1)b_0l_k}
}Y_1 \leq C_g g, \qquad \forall n\geq N_0:=2+\left\lfloor\Big| \log_2
\left(
\frac{C_0 Y_1}{C_g g} \right)/b_0 l_k\Big|\right\rfloor.
\]
Therefore one derives from \eqref{Fi1}
\[
Y_n\leq 2C_g g, \qquad \forall n\geq N_0.
\]
This establishes the desired result.  The proof is complete.

\end{proof}

{\begin{rem} The requirement of $\mathbf{u^0}, \mathbf{u^1}\in \mathbf{V}$ is due to the explicit discretization of $N(\mathbf{u})$ in \eqref{DDSe0} and \eqref{DSAVe0}. This can be achieved by initializing the scheme with energy stable algorithms that allow $\mathbf{H}$ initial data, such as the backward Euler scheme, and utilizing the mr-SAV-BDF2 scheme starting from $n=3$. However, such an approach does not fit the dynamical system approach needed for the convergence of long-time dynamics although it is fine from the purely numerical point of view. The regularity restriction can be relaxed for some systems, as we shall illustrate below for the continuously stratified quasi-geostrophic equation, by utilizing a fractional Sobolev space as the phase space for the numerical scheme. This space is smaller/smoother than the original phase space, but larger than $\mathbf{V}$ so that compactness of the solution in the phase space can still be derived with uniform bound in $\mathbf{V}$. This approach is compatible with the dynamical system framework we adopt. 
\end{rem}}

\section{Uniform large time $H^1$   estimate }\label{sta-H1}

To fix ideas, we perform the analysis on the mean-reverting SAV-BDF2 scheme  for the damped-driven continuously stratified QG model. Throughout this section, one assumes that $f$ is time-independent  so that the autonomous system could reach a statistical equilibrium. { We first introduce the fractional Sobolev function space $H^\alpha_p, \frac14 \le \alpha <1$ for the phase space such that the scheme  is well-posed in the phase space. Then we establish  a large-time $H^1$ bound on the numerical solution with the bound independent of initial data from the phase space and time steps. This estimate is 
 critical to the convergence of the long-term behavior of the system.}

Recall the 3D QG system from \eqref{SVe}--\eqref{BC}.  
Since the anisotropic viscosity in $D(\omega)$ does not affect the analysis in any essential way, for notational simplicity, we replace $D(\omega)$ with $\nu \Delta \omega$ throughout this section.  One introduces the following  expanded continuously stratified QG system
\begin{align}
	&\frac{\partial \omega}{\partial t}-\nu \Delta \omega  +q(t) \big[\nabla^{\perp}_H \psi \cdot \nabla_H \omega+\beta \frac{\partial \psi}{\partial x}\big]=f, \quad \omega =-\big(\Delta_H \psi+F^2 \frac{\partial^2 \psi}{\partial z^2}\big), \label{SVe1}\\
	&\frac{\partial q}{\partial t}+\gamma q-\big<\nabla^{\perp}_H {\psi} \cdot \nabla_H {\omega} +\beta \frac{\partial \psi}{\partial x},\omega\big>=\gamma.  \label{SVe2}
\end{align}
The mr-SAV-BDF2 scheme is  recalled as follows
\begin{align}
	&\delta_t \omega^{n+1}- \nu \Delta \omega^{n+1}+ q^{n+1} \Big[\nabla^{\perp}_H \overline{\psi}^{n+1} \cdot \nabla_H \overline{\omega}^{n+1}+\beta \frac{\partial \overline{\psi}^{n+1}}{\partial x}\Big]=f^{n+1}, \label{HDS-SVe1} \\
	& \omega^{n+1}=-\Big(\Delta_H \psi^{n+1}+F^2 \frac{\partial^2 \psi^{n+1}}{\partial z^2}\Big), \label{HS-SVe12}\\
	&\delta_t q^{n+1}+\gamma q^{n+1}- \Big<\nabla^{\perp}_H \overline{\psi}^{n+1} \cdot \nabla_H \overline{\omega}^{n+1} +\beta \frac{\partial \overline{\psi}^{n+1}}{\partial x},\omega^{n+1}\Big> =\gamma. \label{HDS-SVe2}
\end{align}

First, we address the well-posedness of the mr-SAV-BDF2 scheme in a suitably chosen phase space which is more regular than the phase space of the continuous model.
\begin{lemma}\label{Well-po}
	For $f \in L^2(\Omega)$ and any $\omega^{n-1}, \omega^{n} \in H^\alpha_p$, $\alpha \in [\frac{1}{4}, 1), n\geq 1$, the scheme \eqref{HDS-SVe1}--\eqref{HDS-SVe2} admits a unique solution $\omega^{n+1} \in H^1_{p0}(\Omega)$.
\end{lemma}
\begin{proof}
	Since $\omega^{n-1}, \omega^{n} \in H^\alpha_p$, the elliptic regularity theory implies $\overline{\psi}^{n+1} \in H^{2+\alpha}_p$. Hence  for $\alpha \in [\frac{1}{4}, 1)$, Sobolev embedding implies that $\nabla^{\perp}_H \overline{\psi}^{n+1} \cdot \nabla_H \overline{\omega}^{n+1} = \nabla_H\cdot\left(\nabla^{\perp}_H \overline{\psi}^{n+1} \overline{\omega}^{n+1}\right) \in H^{-1}_{p}$, the dual space of $H^1_{p0}(\Omega)$. Therefore, Eq. \eqref{HDS-SVe2} is well-defined. The existence and uniqueness of a solution $\omega^{n+1} \in H^1_{p0}(\Omega)$ follow easily from  the linear superposition Eqs. \eqref{Decomp}--\eqref{Dec-sol} and standard elliptic theory.  This completes the proof.
\end{proof}

The large time  estimate in the $H^1$ norm of the numerical solution to \eqref{HDS-SVe1}--\eqref{HDS-SVe2} is established in the following theorem.

{ 
\begin{theorem}\label{main-co}
	Let \(f\in L^2(\Omega)\) be time-independent,  and
\[
        (\omega^i,q^i)\in H_p^\alpha(\Omega)\times \mathbb R,
        \qquad i=0,1,
        \qquad \frac14 \le \alpha <1.
\] 
Then there exists a constant  \(M_g>0\) independent of the initial data and $k$, and an integer $N_1(k)$ such that
	\begin{align}\label{gra-es}
		|| \nabla V^{n}||_G^2 \leq R_1^2, \quad \forall n\geq N_1,
	\end{align}	
    where  
\(    
R_1^2:=M_g \left[        R_0^{18}+\|f\|^2+1\right]. 
\)
\end{theorem}
}

\begin{proof} { The arguments below  can be made rigorous via a standard Galerkin approximation.}  We assume $n\geq N_0+1$ such that the estimate \eqref{zero-es} also holds  for $n-1$.

\noindent {  Step 1.  One-step energy estimate.  }  	
	Testing Eq. \eqref{HDS-SVe1} by $-k \Delta \omega^{n+1}$, performing integration by parts, invoking the BDF2-G identity \eqref{G-identity} and H\"older's inequality, one obtains
	\begin{equation}\label{H1es-1}
		\begin{split}
		&|| \nabla V^{n+1}||_G^2-|| \nabla V^{n}||_G^2+\frac{1}{4}||\nabla (\omega^{n+1}-2\omega^n+\omega^{n-1})||^2 +\nu k
		|| \Delta \omega^{n+1}||^2\\
		&\leq \frac{\nu k}{4}  ||\Delta  \omega^{n+1}||^2
		+\frac{C k}{\nu} (||f||^2+{ \beta^2} |q^{n+1}|^2||\overline{\omega}^{n+1}||^2)\\
		&\quad+k \left|q^{n+1} \int_\Omega \nabla^{\perp}_H \overline{\psi}^{n+1} \cdot \nabla_H \overline{\omega}^{n+1} \Delta \omega^{n+1}\, d\mathbf{x}\right|,
		\end{split}
	\end{equation}
    { where the estimate on $||\frac{\partial \overline{\psi}^{n+1}}{\partial x}||$  follows from Eq. \eqref{HS-SVe12}. }

    In light of Eq. \eqref{HS-SVe12}, and by H\"older's inequality,   one controls the trilinear term as follows
	\begin{equation}\label{H1es-2}
		\begin{split}
				&\left|q^{n+1} \int_\Omega \nabla^{\perp}_H \overline{\psi}^{n+1} \cdot \nabla_H\overline{\omega}^{n+1} \Delta \omega^{n+1}\, d\mathbf{x}\right|\\
				&\leq |q^{n+1}| \cdot ||\Delta \omega^{n+1}||\cdot ||\nabla^{\perp}_H \overline{\psi}^{n+1}||_{L^6} || \nabla_H \overline{\omega}^{n+1} ||_{L^3} \\
		&\leq C |q^{n+1}| \cdot  ||\Delta \omega^{n+1}|| (||\omega^n||+||\omega^{n-1}||)||\nabla (2\omega^n-\omega^{n-1}) ||^{\frac{1}{2}} ||\Delta (2\omega^n-\omega^{n-1})||^{\frac{1}{2}} \nonumber \\
		&\leq C |q^{n+1}| \cdot ||\Delta \omega^{n+1}|| (||\omega^n||+||\omega^{n-1}||)||(2\omega^n-\omega^{n-1}) ||^{\frac{1}{4}} ||\Delta (2\omega^n-\omega^{n-1})||^{\frac{3}{4}} \nonumber \\
        & \leq C |q^{n+1}| \cdot ||\Delta \omega^{n+1}|| (||\omega^n||+||\omega^{n-1}||)^{\frac{5}{4}} (||\Delta \omega^n||+||\Delta \omega^{n-1}||)^{\frac{3}{4}} \nonumber \\
        & \leq C ||\Delta \omega^{n+1}|| \cdot |q^{n+1}| (||\omega^n||+||\omega^{n-1}||)^{\frac{5}{4}} (||\Delta \omega^n||^2+||\Delta \omega^{n-1}||^2)^{\frac{3}{8}} \nonumber \\
		&\leq \frac{\nu }{4} ||\Delta \omega^{n+1}||^2+ \frac{C}{\nu^7} |q^{n+1}|^8 (||\omega^n||+|| \omega^{n-1}||)^{10}+ \frac{\nu}{32} ||\Delta \omega^n||^2+\frac{\nu}{32} ||\Delta \omega^{n-1}||^2.
		\end{split}
	\end{equation}
	{   In the above derivation, the second inequality follows from the embedding $H^1\hookrightarrow L^6$, the elliptic regularity theory as well as the interpolation of $L^3$ between $L^2$ and $H^1$ \cite[pp.139, Theorem 5.8]{AdFo2003}; the third inequality is a result of  the interpolation of $H^1$ between $L^2$ and $H^2$; and the last step follows from Young's inequality 
    $$a\cdot r\cdot c\leq \frac{a^2}{2}+\frac{r^8}{8}+\frac{3c^{8/3}}{8},$$
    with $a=\sqrt{\frac{{\nu}}{{2}}} ||\Delta \omega^{n+1}||$, $c=\big(\frac{\nu}{12}\big)^\frac{3}{8} (||\Delta \omega^n||^2+||\Delta \omega^{n-1}||^2)^{\frac{3}{8}}$.
    }

	Substituting the estimate of the nonlinear term into \eqref{H1es-1} and utilizing  the $L^2$ estimate \eqref{zero-es},  one obtains
	\begin{equation}\label{H1es-3}
		\begin{split}
		&|| \nabla V^{n+1}||_G^2-|| \nabla V^{n}||_G^2 +\frac{\nu k}{2}
		|| \Delta \omega^{n+1}||^2 \leq  \frac{k\nu}{32} (||\Delta \omega^n||^2+||\Delta \omega^{n-1}||^2) + k g,
		\end{split}
	\end{equation}
    where 
    $$g:=\frac{C}{ \nu^7}\left(R_0^{18}+1\right)+\frac{C}{\nu}\|f\|^2.$$
\noindent {  Step 2.  Application of Lemma \ref{Lem-key}. }
Define
\[Y_n:=|| \nabla V^{n}||_G^2,\qquad r_n:=\|\Delta\omega^n\|^2,  \qquad
\mu:=\frac{\nu}{2}, \qquad \epsilon:=\frac{1}{16}.
\]
Then \eqref{H1es-3} is written as
\[ Y_{n+1}-Y_n+\mu k r_{n+1}
\le
\varepsilon\mu k(r_n+r_{n-1})+kg,\]
and the condition \eqref{Bas_as} is met with $c_d:=c_0 c_l$,  in light of $||\Delta \omega^n||^2 \ge c_0 ||\nabla \omega^n||^2$ and the norm equivalence \eqref{equi1}. 

An application of Lemma \ref{Lem-key} with $M=N_0+4$ yields
\begin{align}
    \label{Bas_reH1}
    Y_{n+1}
\le
\frac{\mathcal E_M
}{
(1+b_1 l_k)^{\,n-N_0-3}
}
+K_g g, \quad n\geq N_0+4,
\end{align}
where
\[\mathcal E_M:=\max\left\{1,\frac{96}{\nu c_0c_l}\right\} \Big[
Y_M+\frac{5\nu}{32}r_M
+\frac{\nu}{32}  r_{M-1}
\Big], \quad 
b_1
:=\min\left\{1,\frac{\nu c_0 c_l}{48}\right\},  
\]
and \[K_g:=\max\left\{1,\frac{96}{\nu c_0c_l}\right\} \frac{1}{b_1}.\]

\noindent { Step 3. Estimate of \(\mathcal E_M\).}	The elementary $L^2$ energy estimate  \eqref{bas-l2} and the uniform estimate \eqref{zero-es} imply
	\begin{align*}
		||\nabla \omega^n||^2 \leq \frac{C}{\nu k}R_0^2+ \frac{C}{\nu }R_0^2, \qquad n\ge N_0+1.
	\end{align*}
	It follows from the elliptic regularity theory of Eq. \eqref{HS-SVe12} that 
	\begin{align*}
		&||\nabla^{\perp}_H \overline{\psi}^{n+1} \cdot \nabla_H \overline{\omega}^{n+1}||\leq ||\nabla^{\perp}_H \overline{\psi}^{n+1}||_{L^\infty} ||\nabla_H \overline{\omega}^{n+1}|| \leq C||\overline{\psi}^{n+1}||_{H^3} ||\nabla \overline{\omega}^{n+1}|| \\
		&\leq C ||\nabla \overline{\omega}^{n+1}||^2 \\
		&\leq \frac{C}{\nu k}R_0^2+ \frac{C}{\nu }R_0^2, \qquad n\ge N_0+2.
	\end{align*}
	The elliptic regularity estimate of Eq. \eqref{HDS-SVe1} then gives 
	\begin{align*}
		|| \omega^{n+1}||_{H^2}^2 &\leq \frac{C(\nu, \gamma)}{ k^2}(R_0^6+1)+C(\nu, \gamma) \big(R_0^6 +||f||^2\big), \qquad  n\geq N_0+2.
	\end{align*}
    Since $M=N_0+4$, it follows
    \begin{align}
        \label{E5-es}
        \mathcal E_M &\leq \max\left\{1,\frac{96}{\nu c_0c_l}\right\} \Big[|| \nabla V^{M}||_G^2
+\frac{5\nu}{32}\|\Delta\omega^M\|^2+\frac{\nu}{32}\|\Delta\omega^{M-1}\|^2\Big]\\
&    \leq \frac{C(\nu, \gamma)}{ k^2}\big(R_0^6+1\big)  + C(\nu, \gamma) \big(R_0^6 +||f||^2\big):=A_k. \nonumber 
    \end{align}

\noindent { Step 4. Conclusion of the final estimate.} 
Substituting the estimate \eqref{E5-es} into \eqref{Bas_reH1} gives 
	\begin{align*}
		|| \nabla V^{n+1}||_G^2
        &\leq \frac{A_k}{(1+b_1  l_k)^{n-N_0-3}} +K_g g, \quad n\geq N_0+4.
	\end{align*}
    { 
    }
    It follows that
    	\begin{align*}
		|| \nabla V^{n}||_G^2\leq 2K_g g, \quad n\geq N_1:=N_0+5+\left\lfloor\Big| \log_2
\left(
\frac{A_k}{K_g g} \right)/(b_1 l_k)\Big|\right\rfloor.
	\end{align*} 
Thus the estimate \eqref{gra-es} is proven. The proof is complete.
\end{proof}

\begin{rem}
	The fully periodic case can be handled in exactly the same way after we separate the average of the vorticity field. Notice that the mean of the vorticity is an invariant of the dynamics assuming the external forcing is mean-zero. The average of the vorticity contains a linear growth term if the average of the external forcing term is non-zero. Once we focus on the fluctuation (mean-zero) part, we can apply the Poincar\'e-Wirtinger inequality instead of the classical Poincar\'e inequality. 
\end{rem}

\section{Asymptotic consistency}\label{Asym-sim}
The asymptotic consistency refers to the closeness of $\omega^{n+1}$ and $\omega^n$ as well as that of $q^n$ and $1$ for large $n$. This is a foundational piece for convergence of the global attractor and invariant measure  of the mr-SAV-BDF2 scheme to those of the continuous model { since the natural dynamical system formulation of the scheme is on the product phase space $(H\times\mathbb{R})^2$  that contains $(\omega^n,\omega^{n+1}, q^n, q^{n+1})$ while the original model is formulated in terms of $\omega$ only (with $q(t)\equiv 1$).}

{ 
\subsection{Consistency of neighboring points}
We first establish the asymptotic neighboring consistency in the sense that the difference between $(\omega^{n+1}, q^{n+1})$ and $(\omega^{n}, q^{n})$  is small for large $n$. More specifically, we have the following results.
}
    { 
\begin{lemma}[Asymptotic consistency]\label{asym-con}
 For any $k>0$, there exists an integer $N_a$ dependent on $k$ and the initial data such that  for all $n\geq N_a$,
 the solution to the scheme \eqref{HDS-SVe1}--\eqref{HDS-SVe2} satisfies 
	\begin{align}\label{asym-con-eq1}
		||\omega^{n+1}-\omega^n||+k^{{-\frac{1}{2}}}\big(||\omega^{n+1}-\omega^n||_{H^{-1}_p}+|q^{n+1}-q^n|\big) \leq C k^{\frac{1}{2}},
	\end{align}
    where the constant $C$ is independent of $k$ and the initial condition.
\end{lemma}
}
\begin{proof}
One takes advantage of the  compression effect of the BDF2 approximation of the derivative on the time difference. 
	One writes Eq. \eqref{HDS-SVe1} as
	\begin{align*}
		3(\omega^{n+1}-\omega^{n})=( \omega^{n}-\omega^{n-1})+2k \Big(\nu \Delta \omega^{n+1}-q^{n+1} \big[\nabla^{\perp}_H \overline{\psi}^{n+1} \cdot \nabla_H \overline{\omega}^{n+1}+\beta\overline{\psi}^{n+1}_x\big]+f\Big).
	\end{align*}
    To simplify the computation, throughout the proof, { one assumes that $n\geq N:=N_1+1$  such that the estimates \eqref{zero-es} and \eqref{gra-es} hold.} 
Then
	\begin{align*}
		&3||\omega^{n+1}-\omega^{n}||_{H^{-1}_p} \nonumber \\
		&\leq ||\omega^{n}-\omega^{n-1}||_{H^{-1}_p}+Ck ||\nu \Delta \omega^{n+1}-q^{n+1} \big[\nabla^{\perp}_H \overline{\psi}^{n+1} \cdot \nabla_H \overline{\omega}^{n+1}+\beta\overline{\psi}^{n+1}_x\big]+f||_{H^{-1}_p} \nonumber\\
		&\leq  ||\omega^{n}-\omega^{n-1}||_{H^{-1}_p}+C k \Big(||\nabla \omega^{n+1}||+|q^{n+1}| \cdot || \overline{\omega}^{n+1}||\cdot  ||\nabla \overline{\omega}^{n+1}|| +||\overline{\omega}^{n+1}||+||f|| \Big) \nonumber\\
		&\leq  ||\omega^{n}-\omega^{n-1}||_{H^{-1}_p}+Ck,
	\end{align*}
    {  where the constant $C$ changes from line to line but is independent of $k$ and the initial condition.}
An iteration of the above  inequality gives
	\begin{equation}\label{con-1}
		\begin{split}
		||\omega^{n+1}-\omega^{n}||_{H^{-1}_p}  
		&\leq \frac{1}{3^{n-N}} ||\omega^{N+1}-\omega^{N}||_{H^{-1}_p} +\frac{3}{2}C k   \\
		&\leq \frac{C R_0}{3^{n-N}} +\frac{3}{2}C k,
		\end{split}
	\end{equation}
where the last step follows  from Theorem \ref{theo1}. Likewise,
\begin{equation}\label{con-1q}
		|q^{n+1}-q^n|  \leq \frac{C R_0}{3^{n-N}} +\frac{3}{2}C k.
	\end{equation}
    Hence for all \(n\geq N+\left\lfloor\big|\ln{\frac{C R_0}{k}}/{\ln 3}\big|\right\rfloor\),
\[||\omega^{n+1}-\omega^{n}||_{H^{-1}_p}+|q^{n+1}-q^n|\leq Ck.\] 
	Through interpolation between $H^{-1}_p$ and $H^1_{p0}$, one obtains
	\begin{align*}
		||\omega^{n+1}-\omega^{n}||^2\leq ||\omega^{n+1}-\omega^{n}||_{H^{-1}_p}||\omega^{n+1}-\omega^{n}||_{H^{1}} \leq C k.
	\end{align*}
This completes the proof.  
\end{proof}

\begin{rem}
	The $\|\omega^{n+1}-\omega^{n}\|_{H^{-1}}$ and $|q^{n+1}-q^n|$ estimates are optimal since they approximate the true solution at two neighboring time levels with a distance of $k$.
\end{rem}

{ 
\subsection{Long time consistency of the scalar auxiliary variable}
Next, we show that the scalar auxiliary variable $q^n$ is close to the desired value of 1 for all large time $n$. 
This is a pleasant surprise since we usually anticipate the accumulation of small errors. 
The closeness of the auxiliary variable to 1 is one of the foundational pieces of the convergence of the long-time dynamics.
}
    { 
\begin{lemma}\label{asym-con-1} For any $k>0$, 
there exists an integer $N_3$ dependent on $k$ and the initial data such that for any $n\geq N_3$,
	\begin{align}\label{DS-SVe2n3Th}
		|q^n-1| \leq C k^{\frac{1}{2}},
	\end{align} 
where the constant $C$ is independent of $k$ and the initial condition.
\end{lemma}
}
\begin{proof}

  Define $a^n=q^n-1$.  Since the vorticity  is in  $H^1$, the duality in Eq. \eqref{HDS-SVe2} can be interpreted as integration and hence it can be reformulated as, with the help of the skew symmetry of the advection term and the $\beta$ term, 
	\begin{align}\label{DS-SVe2n}
		&\delta_t a^{n+1}+\gamma a^{n+1}=\int_\Omega  \big[\nabla^{\perp}_H \overline{\psi}^{n+1} \cdot \nabla_H \overline{\omega}^{n+1}+ \beta\overline{\psi}^{n+1}_x\big]\omega^{n+1}\, d{\mathbf{x}} \nonumber \\
		&=\int_\Omega  \nabla^{\perp}_H \overline{\psi}^{n+1} \cdot \nabla_H\big[\overline{\omega}^{n+1} -\omega^{n+1}\big]\omega^{n+1}+\beta[\overline{\psi}^{n+1}_x-{\psi}^{n+1}_x]\omega^{n+1}\, d{\mathbf{x}}.
	\end{align}
	 Assuming that     { $n\geq N_a+1$} as in Lemma  \ref{asym-con}, one derives  from Theorem \ref{main-co} and Lemma \ref{asym-con} that 
	\begin{align*}
		&\Big| \int_\Omega  \nabla^{\perp}_H \overline{\psi}^{n+1} \cdot \nabla_H \big[\overline{\omega}^{n+1} -\omega^{n+1}\big]\omega^{n+1}\, d{\mathbf{x}}\Big|\\
        &= \Big| \int_\Omega  \nabla^{\perp}_H \overline{\psi}^{n+1} \cdot \nabla_H \omega^{n+1} \big[\overline{\omega}^{n+1} -\omega^{n+1}\big]\, d{\mathbf{x}}\Big|\\
        &\leq ||\nabla^{\perp}_H \overline{\psi}^{n+1}||_{\infty} \cdot ||\nabla \omega^{n+1}|| \cdot ||\overline{\omega}^{n+1} -\omega^{n+1}||\\
        &\leq C ||\nabla \overline{\omega}^{n+1}||\cdot ||\nabla \omega^{n+1}|| \cdot ||\overline{\omega}^{n+1} -\omega^{n+1}||\\
        &\leq C ||\overline{\omega}^{n+1} -\omega^{n+1}||\\
        &\leq C k^\frac12.
	\end{align*}
	Therefore  multiplying Eq. \eqref{DS-SVe2n} by $ka^{n+1}$ and utilizing Young's inequality one obtains
	\begin{align*}
		|A^{n+1}|^2_G+\frac{\gamma}{2} k |a^{n+1}|^2 
		&\leq |A^{n}|^2_G+\frac{C k}{\gamma} ||\overline{\omega}^{n+1} -\omega^{n+1}||^2  \nonumber \\
		&\leq   |A^{n}|^2_G  
		 +C k^2,
	\end{align*}
	with $A^n=[a^n, a^{n-1}]^T$. Applying Lemma \ref{Lem-key} to the above inequality with 
\[ M=N_a+1, \quad Y_n:=|A^{n}|^2_G, \quad r_n:=|a^{n}|^2, \quad \mu=\frac{\gamma}{2},  \quad g=Ck, \quad \epsilon=0,\] one has for $n\geq N_a+1$,
\begin{align*}
   |A^{n+1}|^2_G&\leq \frac{C}{\rho_k^{\,n-N_a}}\Big(\|A^{N_a+1}\|_G^2+\gamma/8|a^{N_a+1}|^2\Big)
        +       C k, \\
        &\le \frac{C(R_0^2+1)}{\rho_k^{\,n-N_a}}+Ck,
\end{align*} 
where \(b_3:=\min\left\{1,\frac{\gamma c_l}{48}\right\}\) and
\( \rho_k:=1+b_3 l_k,\) and the last step is a consequence of the $L^2$ estimate \eqref{zero-es}.

Therefore for all $n\geq N_a+2+\left\lfloor\big|\log_2{\big[C(R_0^2+1)/k\big]} /(b_3 l_k)\big|\right\rfloor:=N_3$,
\[|q^n-1|^2   \leq C k.\]
This completes the proof.
\end{proof}

{ 
\subsection{mr-SAV-BDF2 as a perturbation of IMEX-BDF2}
With the uniform estimate on $q^n-1$, we can show that our mr-SAV-BDF2 scheme is a small perturbation of the classical IMEX-BDF2 scheme when one treats the dissipation term implicitly and the nonlinear term explicitly.}
	The general  mr-SAV-BDF2 scheme \eqref{DDSe0} can be re-written as
	\begin{align*}
		\delta_t \ub^{n+1} +A  \ub^{n+1}+ N\big(\overline{\ub}^{n+1}\big) =(1-q^{n+1}) N\big(\overline{\ub}^{n+1}\big)+\mathbf{F}^{n+1}, 
	\end{align*}
	where the left-hand side is the classical     { IMEX-BDF2 scheme with  explicit discretization of the nonlinear term}. For the continuously stratified QG equation, 
	the asymptotic estimate \eqref{DS-SVe2n3Th} and the large-time uniform $H^1$ estimate on the numerical solution imply that for  $n$ sufficiently large
	\begin{align*}
	&	||(1-q^{n+1}) N\big(\overline{\ub}^{n+1}\big)||_{L^2} \\
	&\leq |1-q^{n+1}| \left(\|\nabla^{\perp}_H \overline{\psi}^{n+1} \cdot\nabla_H\overline{\omega}^{n+1}\|_{L^2} +|\beta|\|\frac{\partial\overline{\psi}^{n+1}}{\partial x}\|_{L^2}\right)\\
	&\le C_1 k^{\frac{1}{2}}. 
	\end{align*}
   This suggests that in the asymptotic regime the mr-SAV-BDF2 scheme 
	is close to the standard BDF2 scheme with fully explicit discretization of the nonlinear term while preserving the uniform-in-time estimates regardless of the time step size.  {Such an asymptotic estimate as \eqref{DS-SVe2n3Th} is not available, and perhaps impossible as suggested by our numerical results presented below, to other SAV type schemes.}

\section{Convergence of attractors and invariant measures}\label{long-time}
For dissipative dynamical systems with the global attractor, invariant measures are supported on the global attractor. The global attractor and invariant measure completely characterize the long-time dynamics of the underlying dynamical system. In this section we show that the mr-SAV-BDF2 scheme captures the long-time dynamics of the continuously stratified QG model. We first recall the notions of the global attractor and invariant measures.
\subsection{Preliminaries}
Recall that the model \eqref{SVe} generates a dissipative continuous dynamical system $\{\mathbb{S}(t), t\geq 0\}$ on $H=L^2(\Omega)$  which possesses the global attractor $\mathcal{A} \subset H$ \cite{Wang1992, Temam1997}.

{ 
\begin{definition}[Global Attractor]\label{def-att}
    A compact set $\mathcal{A} \subset H$ is the global attractor of $\mathbb{S}(t)$ if the following conditions are satisfied:
    \begin{enumerate}[(i).]
        \item $\mathcal{A}$ is invariant under the flow, i.e. $$\mathbb{S}(t)\mathcal{A}=\mathcal{A}, \quad \forall t\geq 0;$$
        \item $\mathcal{A}$ attracts all bounded sets in $H$, i.e., for every bounded set $B\subset H$,
        $$ \lim_{t\rightarrow \infty} \delta_H \big(\mathbb{S}(t) B, \mathcal{A}\big)=0,$$
        where  $\delta_H(A, B)=\sup_{a\in A} \inf_{b\in B}||a-b||_H$ is   the Hausdorff set semi-distance in $H$.
    \end{enumerate}
\end{definition}
}
Similarly, the expanded system \eqref{SVe1}--\eqref{SVe2} induces a continuous semi-group $\{\mathbb{S}_q(t), t\geq 0\}$ 
on the product space $L^2\times \mathbb{R}$. Since $q(t)\rightarrow 1, t\rightarrow +\infty$ for any initial condition $q_0$ and $q(t)\equiv1$ if $q_0=1$, the global attractor $\mathcal{A}_q$ of $\mathbb{S}_q(t)$ must be equal to $\mathcal{A}\times \{1\}$.

 Next, we introduce the notion of stationary statistical solution (invariant measure) for the continuously stratified QG equation. Cf. \cite[pp.183]{FMRT2001} for the corresponding definition  to the 2D Navier-Stokes equations, and \cite{Wang2012}  for the 2D NSE in the vorticity--stream function formulation.{ 
\begin{definition}[Stationary Statistical Solution]\label{Sta-def}
 A Borel probability measure $\mu$ on $H$ is  a stationary statistical solution of the continuously stratified QG equation \eqref{SVe}, if
     \begin{align*}
         &\int_{H} ||\nabla \omega||^2 d\mu <\infty, \\
         &\int_{H} 
		\Big\langle-f-\nu \Delta \omega + \nabla^{\perp}_H \psi \cdot \nabla_H \omega+\beta\psi_x, \Phi^\prime(\omega)\Big\rangle\, d \mu=0, \\
        &\int_{H} \nu  ||\nabla  \omega||^2-\big(f, \omega\big)\, d\mu \leq 0,
     \end{align*}
     where $\langle, \rangle$ denotes the duality pairing between $H^{-1}_p$ and $H^1_{p0}$.
{  Here $\Phi$ is an arbitrary cylindrical test function in the form of}
\[\Phi(\omega):=\phi\big((\omega, g_1), \cdots, (\omega, g_m)\big)\] 
where $\phi$ is a smooth function of compact support on $\mathbb{R}^m$,   and $g_i, i=1, 2,\ldots m$ are the first $m$ eigenfunctions of the Laplace operator with periodic boundary condition in the horizontal direction and homogeneous Dirichlet boundary condition at $z=0,1$. 
\end{definition}
\begin{definition}[Invariant Measure]
   A Borel probability measure $\mu$ on $H$ is called an invariant measure of the dynamical system generated by the semi-group $\mathbb{S}(t), t\geq0$ if
   $$
   \int_H \Phi (\mathbf{u}) d\mu =\int_H \Phi\big(\mathbb{S}(t) \mathbf{u}\big) d\mu, \quad \forall t\geq0,
   $$
   for all bounded continuous test functionals $\Phi$.
\end{definition}}
For the continuously stratified QG model, it can be shown that the invariant measures are equivalent to stationary statistical solutions, cf. \cite[Chap. IV, Appendix B.1]{FMRT2001} for a rather technical proof in the case of 2D Navier-Stokes equations.
In the Definition \ref{Sta-def}, the first condition implies that the invariant measure is supported in the finer space $H^1_{p0}$; the second condition is the weak formulation of the invariance of the measure under the flow; and the third condition is a statistical version of the energy inequality.

{Let $\mathcal{IM}$ and $\mathcal{IM}_q$ be the sets of all invariant measures of the damped-driven continuously stratified quasi-geostrophic equations \eqref{SVe} and the expanded PDE system \eqref{SVe1}--\eqref{SVe2}, respectively. 
 We recall that the continuous model when viewed as a dynamical system on $H^\alpha,  1/4\le \alpha <1$ has the same global attractor and invariant measures as with $L^2$ as the phase space, since  all invariant measures are supported on the global attractor \cite{Wang2009c} { which is a subset of $H^1$}.

\subsection{Dynamical system formulation of the scheme and projections}
{ In order to apply our Lax-type result on convergence of invariant measures and global attractors, we formulate our two-step mr-SAV-BDF2 scheme as a discrete dynamical system.}

Thanks to Lemma \ref{Well-po}, the scheme \eqref{HDS-SVe1}--\eqref{HDS-SVe2} can be written as a dynamical system on the product space $(H^\alpha_p \times \mathbb{R})^2, \alpha \in [\frac{1}{4}, 1)$:
\begin{align}\label{dyna}
	\mathbb{S}_{q,k} \begin{bmatrix}
		(\omega^{n}, q^{n})\\
		(\omega^{n-1}, q^{n-1})
	\end{bmatrix}=\begin{bmatrix}
		(\omega^{n+1}, q^{n+1})\\
		(\omega^{n}, q^{n})
	\end{bmatrix}.
\end{align}
{ It is clear that $\mathbb{S}_{q,k}$ is continuous on the product space $(H^\alpha_p \times \mathbb{R})^2, \alpha \in [\frac{1}{4}, 1)$.}
The following corollary follows directly from the estimate \eqref{zero-es} in Theorem \ref{theo1},  the  estimate  \eqref{gra-es} in Theorem \ref{main-co} and the compact embedding $H^1 \hookrightarrow H^\alpha_p, \alpha \in [\frac{1}{4}, 1)$.
\begin{cor}\label{Uni-Cor}
    For any $k>0$,  the dynamical system $\mathbb{S}_{q,k}$ on the product space $(H^\alpha_p \times \mathbb{R})^2$, $\alpha \in [\frac{1}{4}, 1)$ possesses the global attractor $\mathcal{A}_{q, k}$ such that
	\begin{equation}\label{att-boun}
		\begin{split}
	&|| w_1||_{H^1}^2+|| w_2||_{H^1}^2+|q_1|^2+|q_2|^2 \\
	&\leq C(R_0^2+R_1^2), \quad \forall \left[(w_1, q_1), (w_2, q_2)\right]^T \in \mathcal{A}_{q, k}.
	\end{split}
	\end{equation} 
\end{cor}
The set of all invariant measures for $\mathbb{S}_{q,k}$ on the product space $(H^\alpha_p \times \mathbb{R})^2$ is denoted by $\mathcal{IM}_{q,k}$.  $\mathcal{IM}_{q,k}$ is non-empty thanks to the $H^1$ uniform boundedness of the global attractor $\mathcal{A}_{q,k}$ and by the Bogliubov-Krylov construction of invariant measures \cite{FMRT2001}.

{The dynamical system \eqref{dyna} generated by the scheme \eqref{HDS-SVe1}--\eqref{HDS-SVe2} takes two arbitrary sets of initial data $(\omega^0, q^0)$ and $(\omega^1, q^1)$ and lives on the product space of $(H^\alpha\times\mathbb{R})^2$. On the other hand, the expanded PDE system \eqref{SVe1}--\eqref{SVe2} requires only one initial datum and lives on $H^\alpha\times\mathbb{R}$. 
Therefore, the two components of the dynamical system \eqref{dyna} must be close to each other if the long-time dynamics of the scheme has a chance of approximating those of the original model. The following corollary, which states that 
 the two components of the global attractor $\mathcal{A}_{q, k}$ are indeed close to each other,  is a consequence of the asymptotic consistency established in Lemma \ref{asym-con} and \ref{asym-con-1} as well as the invariance of the global attractor. 
}

\begin{cor}\label{asym-conCor}
 There exists a constant $C$ independent of $k$ such that for any $\left[(\omega_1, q_1), (\omega_2, q_2)\right]^T \in \mathcal{A}_{q, k}$ there holds
	\begin{align}\label{asym-con-eq}
		||\omega_1-\omega_2||+k^{{-\frac{1}{2}}}\big[||\omega_1-\omega_2||_{H^{-1}}+|q_1-q_2|\big]+|q_1-1| \leq C k^{\frac{1}{2}}.
	\end{align}   
\end{cor}

{ 
Notice that the global attractor and the invariant measures of the continuous model live in $H^1$ while the global attractor and the invariant measures of the discrete dynamical system \eqref{dyna} induced by the mr-SAV-BDF2 scheme live on the expanded product space $(H^\alpha_p \times \mathbb{R})^2$, we will need to invoke projection and/or lift  in order to compare the two objects.
Denote $\mathbb{P}_j, j=1, 2$ the projection operator from the product space $(H^\alpha_p \times \mathbb{R})^2$ onto the $j$th coordinate $H^\alpha \times \mathbb{R}$, {  e.g.,
\[\mathbb{P}_1 \begin{bmatrix}{(\omega_1,q_1)}\\{(\omega_2,q_2)}\end{bmatrix}=(\omega_1,q_1). \]
} The induced projection operator $\mathbb{P}_j^\ast$ in the space of Borel probability measures is defined as
\begin{align*}
	\mathbb{P}_j^\ast \mu(B)= \mu(\mathbb{P}_j^{-1} B), \quad \forall B \in \mathcal{B}(H^\alpha_p \times \mathbb{R}),
\end{align*}
where $\mathcal{B}(H^\alpha_p \times \mathbb{R})$ is the Borel $\sigma-$algebra of $H^\alpha_p \times \mathbb{R}$, and $\mu$ is a Borel probability measure on the product space. Likewise,  $\mathbb{P}_\omega$ is understood as the coordinate projection from $H^\alpha_p \times \mathbb{R}$ onto $H^\alpha$, and the induced marginal distribution of the measure in the $\omega$ component is denoted by $\mathbb{P}_\omega^\ast$:{ 
\[\mathbb{P}_\omega\Big((\omega,q)\Big):=\omega, \qquad \mathbb{P}_\omega^\ast \lambda (A):=\lambda (\mathbb{P}_\omega^{-1}A), \quad \forall A\in \mathcal{B}(H^\alpha_p),\]
where $\lambda$ is a Borel probability measure on $H^\alpha_p \times \mathbb{R}$.}
}

\subsection{Finite-time uniform convergence}
{In order to prove the convergence of the long-time behavior of the scheme to those of the underlying model using dynamical system approach via Lax-type criteria \cite{Wang2010a,Wang2016}, we need to show that the solutions to the numerical scheme converge to the corresponding solution of the underlying model uniformly for initial data taken from the union of global attractors of the discrete dynamical system \eqref{dyna}.
The next lemma provides the desired convergence albeit with a sub-optimal rate.   
\begin{lemma}[Uniform finite-time trajectory convergence with $H^1$ initial data]\label{Traj-convL}
	For any $T>0$ there exists a constant $C(R_0, R_1, T, \nu)$ such that $\forall \left[(\omega^1, q^1), (\omega^0, q^0)\right]^T \in \mathcal{A}_{q, k}$,  
	\begin{align}\label{Tra-conv}
		||\omega\big(t^{n}\big)-\omega^{n}||_{H^\alpha} \leq C k^{\frac{1-\alpha}{4}}, \quad  t^n=nk \leq T, \quad \frac{1}{4}\leq \alpha<1,
	\end{align}
 {    where $\omega(t)$ is the solution to the expanded system \eqref{SVe1}--\eqref{SVe2} with the initial condition $(\omega^0, q^0)$.}
\end{lemma}
{  Again, the discrete dynamical system  \eqref{dyna} needs two initial conditions \linebreak $\left[(\omega^1, q^1), (\omega^0, q^0)\right]^T$. Despite the arbitrariness of $(\omega^1, q^1)$, in finite time the discrete numerical solution $\{\omega^n\}$ is expected to converge to the continuous trajectory $\omega(t)$ initialized only with  $(\omega^0, q^0)$, due to the asymptotic consistency \eqref{asym-conCor} on the global attractor.  The convergence rate is not optimal because of the arbitrary choice of initial data which may not satisfy the local truncation error required for second order convergence for two-step schemes. The proof is deferred to the Appendix. }

{ 
This finite-time uniform convergence together with the asymptotic consistency  \eqref{asym-con-eq1}, \eqref{DS-SVe2n3Th} implies that the numerical solution converges to the trajectory of  $\mathbb{S}(t)=\begin{bmatrix}
		S_q(t)\\
		S_q(t)
	\end{bmatrix}$ on $\big(H^\alpha_p\times \mathbb{R}\big)^2$ uniformly  in initial data from the discrete global attractor. More specifically, in light of $q(t)=(q_0-1)e^{-\gamma t}+1$, by the triangle inequality and the interpolation inequality, we have
\begin{align}
    \label{Traj-convL1}
    \left\|\begin{bmatrix}
		\big(\omega^{n+1}, q^{n+1}\big)\\
		\big(\omega^{n}, q^{n}\big)
	\end{bmatrix}-\begin{bmatrix}
		\big(\omega_1 (t^{n}), q_1(t^{n})\big)\\
		\big(\omega_0(t^{n}), q_0(t^{n})\big)
	\end{bmatrix}
    \right\|_{\big(H^\alpha_p\times \mathbb{R}\big)^2} \leq C k^{\frac{1-\alpha}{4}},
\end{align}
where \(\big(\omega_i(t), q_i(t)\big), i=0, 1\) are the solutions to the expanded system initiated by \(\big(\omega^i, q^i\big)\), respectively.

}

\subsection{Attractor convergence}
The attractor convergence is relatively straightforward with the uniform bounds in hand, the local in time uniform convergence, the asymptotic consistency, and the relationship between the global attractor of the expanded system and the original one.} 
\begin{theorem}[Attractor convergence]\label{att-conv}
	Let $\{\mathcal{A}_{q,k}\}$ be the family of global attractors  for the semigroup $\mathbb{S}_{q,k}$ in $(H^\alpha_p\times \mathbb{R})^2$, $\alpha \in [\frac{1}{4}, 1)$,  and recall that $\mathcal{A}_q$ and $\mathcal{A}$ are the global attractors associated with the semigroup $S_q(t)$ via the expanded system \eqref{SVe1}--\eqref{SVe2} and the original damped-driven continuously stratified QG model, respectively. Define $\mathds{1}\mathcal{A}_q = \left\{ \big[(\omega, 1), (\omega, 1)\big]^T \in \mathcal{A}_q \times \mathcal{A}_q \right\}$. Then
	\begin{align}
		&\delta_{\big(H^\alpha\times \mathbb{R}\big)^2}\big(\mathcal{A}_{q,k}, \mathds{1}\mathcal{A}_{q}\big)\longrightarrow 0, \quad k\rightarrow 0, \label{Attr-conv-1} \\
		&\delta_{H^\alpha}\big(\mathbb{P}_\omega \mathbb{P}_j \mathcal{A}_{q,k}, \mathcal{A}\big)\longrightarrow 0, \quad k\rightarrow 0, \quad j=1, 2,\label{Attr-conv-2}
	\end{align}
	where $\delta_{X} (A, B):=\sup_{a\in A} \inf_{b\in B}||a-b||_X$ is the Hausdorff semi-distance of $X$.
\end{theorem}
\begin{proof}
	Since $\mathcal{A}_q \times \mathcal{A}_q$ is the global attractor  for $\mathbb{S}(t)=\begin{bmatrix}
		S_q(t)\\
		S_q(t)
	\end{bmatrix}$ on $\big(H^\alpha_p\times \mathbb{R}\big)^2$ and $\mathcal{A}_{q,k}$ is uniformly bounded by \eqref{att-boun}, one deduces, 
	for any $\varepsilon >0, \exists T(\varepsilon, R_0, R_1) > 0$ such that
	\begin{equation}\label{Att-traj-conv0}
		\delta_{\big(H^\alpha\times \mathbb{R}\big)^2}\left(\mathbb{S}(t) \mathcal{A}_{q,k}, \mathcal{A}_q \times \mathcal{A}_q \right) \leq \varepsilon, \quad \forall t \geq T(\varepsilon, R_0, R_1).
	\end{equation} 
	Let $n_T$ be the least integer such that $n_T k \geq T(\varepsilon, R_0, R_1)$.  For any $\big[{V}^1, {V}^0\big]^T \in \mathcal{A}_{q, k}$, there exists $\big[\widetilde{{V}}^1, \widetilde {{V}}^0\big]^T\in \mathcal{A}_{q, k}$ such that 
	$
	\mathbb{S}_{q,k}^{n_T} \begin{bmatrix}
		\widetilde{{V}}^1\\
		\widetilde {{V}}^0
	\end{bmatrix}=\begin{bmatrix}
		{V}^1\\
		{V}^0
	\end{bmatrix}
	$, thanks to the invariance of $\mathcal{A}_{q, k}$ under $\mathbb{S}_{q,k}$. The trajectory convergence \eqref{Traj-convL1} then implies
	\begin{align}
		\label{Att-traj-conv1}
		dist\Big(\big[{V}^1, {V}^0\big]^T, \mathbb{S}(n_T k)\big[\widetilde{{V}}^1, \widetilde {{V}}^0\big]^T\Big) &=dist\Big(\mathbb{S}_{q,k}^{n_T}\big[\widetilde{{V}}^1, \widetilde {{V}}^0\big]^T, \mathbb{S}(n_T k)\big[\widetilde{{V}}^1, \widetilde {{V}}^0\big]^T\Big) \nonumber \\
		&\leq C k^{\frac{1-\alpha}{4}},  \quad \frac{1}{4}\leq\alpha<1,
	\end{align}
	where $dist$ is the distance in the $\big(H^\alpha_p\times \mathbb{R}\big)^2$ norm.
	It follows from \eqref{Att-traj-conv1} and \eqref{Att-traj-conv0} that
	\begin{align}
		\label{Att-traj-conv2}
		&dist\Big(\big[{V}^1, {V}^0\big]^T, \mathcal{A}_q \times \mathcal{A}_q\Big) \nonumber \\
		 &\leq dist\Big(\big[{V}^1, {V}^0\big]^T, \mathbb{S}(n_T k)\big[\widetilde{{V}}^1, \widetilde {{V}}^0\big]^T \Big)+dist\Big(\mathbb{S}(n_T k)\big[\widetilde{{V}}^1, \widetilde {{V}}^0\big]^T, \mathcal{A}_q \times \mathcal{A}_q\Big) \nonumber \\
		&\leq C k^{\frac{1-\alpha}{4}}+\epsilon,
	\end{align}
	whence
	\begin{align}
		\label{Att-traj-conv3}
		\delta_{\big(H^\alpha_p\times \mathbb{R}\big)^2}\Big(\mathcal{A}_{q,k}, \mathcal{A}_q \times \mathcal{A}_q\Big) \leq C k^{\frac{1-\alpha}{4}}+\epsilon.
	\end{align}
	Sending $k\rightarrow 0$ in \eqref{Att-traj-conv3} and in light of the fact that $\epsilon$ is arbitrary,  one obtains
	\begin{align}
		\label{Att-traj-conv4}
		\lim_{k\rightarrow 0}\delta_{\big(H^\alpha_p\times \mathbb{R}\big)^2}\Big(\mathcal{A}_{q,k}, \mathcal{A}_q \times \mathcal{A}_q\Big)=0.
	\end{align}
	The convergence results \eqref{Attr-conv-1} and \eqref{Attr-conv-2} then follow from the asymptotic consistency estimate \eqref{asym-con-eq} as well as the fact that $\mathcal{A}_q=\mathcal{A}\times\{1\}$. This completes the proof.
\end{proof}

\subsection{Convergence of invariant measures}
The convergence of invariant measures, after taking appropriate projection, is established in the following theorem.
\begin{theorem}[Convergence of invariant measures]\label{im-conv}
	\label{Conv-inva}
	Let $f$ be a time-independent forcing. Let $\{\mu_k \in \mathcal{IM}_{q,k}, 0<k<1\}$ be the collection of invariant measures of the numerical scheme \eqref{HDS-SVe1}--\eqref{HDS-SVe2}. Then any sequence $(\mu_k)_{k\in(0,1)}$ contains a subsequence, still denoted by $(\mu_k)_{k\in(0,1)}$, weakly convergent to a probability measure $\mu_0$ on the product space $(H^\alpha_p \times \mathbb{R})^2$. Furthermore, 
	\begin{align}\label{conv-mea}
		\mathbb{P}_\omega^\ast \mathbb{P}_j^\ast \mu_k \rightharpoonup \mathbb{P}_\omega^\ast \mathbb{P}_j^\ast \mu_0 \in \mathcal{IM}, \quad k\rightarrow 0, \quad j=1, 2.
	\end{align}
\end{theorem}
\begin{proof}
{  We break the proof down into three steps.}

\noindent {  Step 1. Tightness of the invariant measures.}

	Recall that the invariant measures are carried by the global attractor, and the global attractor $\mathcal{A}_{q, k}$ of $\mathbb{S}_{q, k}$ is uniformly bounded with respect to $k$ in $(H^1\times \mathbb{R})^2$, cf. Corollary \ref{Uni-Cor}. By the compact embedding $H^1_{p0}(\Omega) \hookrightarrow H^\alpha_p$, $\alpha \in [\frac{1}{4},1)$, one concludes that the set of invariant measures $\{\mu_k \in \mathcal{IM}_{q,k}, 0<k<1\}$ is tight in the space of probability measures on the product space $(H^\alpha \times \mathbb{R})^2$. Hence for any sequence $(\mu_k)_{k\in(0,1)} \subset \{\mu_k \in \mathcal{IM}_{q,k}, 0<k<1\}$ there exists a subsequence, still denoted by $(\mu_k)_{k\in(0,1)}$, and a Borel probability measure $\mu_0$ on $(H^\alpha_p \times \mathbb{R})^2$ such that $\mu_k \rightharpoonup \mu_0$ weakly as $k\rightarrow 0$. 

    \noindent {  Step 2. Weak stationarity.}
    
	In the following, one proves that the marginal distribution $\mathbb{P}_\omega ^\ast\mathbb{P}_1^\ast \mu_0$ is an invariant measure of the  continuously stratified quasi-geostrophic equations \eqref{SVe}, { in the sense of  Definition \ref{Sta-def}.}   The case of $j=2$ is similar. 

    { 
	Let $\Phi(\omega):=\phi\big((\omega, g_1), \cdots, (\omega, g_m)\big)$ be a cylindrical test function introduced in Definition \ref{Sta-def}. { We note that $\Phi$ is smooth and of finite dimension. }
	It follows that the Fr\'echet derivative is given by
	\begin{align}\label{cylin}
		\Phi^\prime(\omega)=\sum_{i=1}^m \frac{\partial}{\partial z_i}\phi\big((\omega, g_1), \cdots, (\omega, g_m)\big) g_i,
	\end{align}
	where $ \frac{\partial}{\partial z_i}$ is the partial derivative  with respect to the $i$-th component of $\phi$.} 
	
	Denote ${V}^i:=(\omega^i, q^i) \in H^\alpha_p\times \mathbb{R}$. Recall that the semi-group generated by the scheme gives
	\begin{align}\label{dyna2}
		\mathbb{S}_{q,k} \begin{bmatrix}
			{V}^1\\
			{V}^0
		\end{bmatrix}=\begin{bmatrix}
			{V}^2\\
			{V}^1
		\end{bmatrix}.
	\end{align}
	In light of the notion of marginal distribution and weak convergence of $(\mu_k)_{k\in(0,1)}$, one obtains
	\begin{align}
		\label{Weak-form1}
		&\int_{H^\alpha_p} 
		\Big\langle-f-\nu \Delta \omega + \nabla^{\perp}_H \psi \cdot \nabla_H \omega+\beta\psi_x, \Phi^\prime(\omega)\Big\rangle \, d\mathbb{P}^\ast_\omega\mathbb{P}^\ast_1 \mu_0 \nonumber \\
		&=\int_{H^\alpha_p\times \mathbb{R}} 
		\Big\langle-f-\nu \Delta \omega + \nabla^{\perp}_H \psi \cdot \nabla_H \omega+\beta\psi_x, \Phi^\prime(\omega)\Big\rangle \, d\mathbb{P}^\ast_1 \mu_0 \nonumber \\
		&=\int_{(H^\alpha_p\times \mathbb{R})^2} 
		\Big\langle-f-\nu \Delta \omega_1 + \nabla^{\perp}_H \psi_1 \cdot \nabla_H \omega_1+\beta\psi_{1x}, \Phi^\prime(\omega_1)\Big\rangle \, d \mu_0   \left(\begin{bmatrix}
			{V}^1\\
			{V}^0
		\end{bmatrix} \right)\nonumber \\
		&=\lim_{k\rightarrow 0}\int_{(H^\alpha_p\times \mathbb{R})^2} 
		\Big\langle-f-\nu \Delta \omega_1 + \nabla^{\perp}_H \psi_1 \cdot \nabla_H \omega_1+\beta\psi_{1x}, \Phi^\prime(\omega_1)\Big\rangle \, d \mu_k   \left(\begin{bmatrix}
			{V}^1\\
			{V}^0
		\end{bmatrix} \right), \nonumber \\
		&=\lim_{k\rightarrow 0}\int_{(H^\alpha_p\times \mathbb{R})^2} 
		\Big\langle-f-\nu \Delta \omega_1 + q_1[\nabla^{\perp}_H \psi_1 \cdot \nabla_H \omega_1+\beta\psi_{1x}], \Phi^\prime(\omega_1)\Big\rangle \, d \mu_k   \left(\begin{bmatrix}
			{V}^1\\
			{V}^0
		\end{bmatrix} \right),
	\end{align}
	where $\omega =-\big(\Delta_H \psi+F^2 \frac{\partial^2 \psi}{\partial z^2}\big)$, and the asymptotic consistency estimate \eqref{asym-con-eq} is utilized in the last step.  One writes
	\begin{equation}\label{refo1}
		\begin{split}
		&-\nu\Delta \omega_1  + q_1[\nabla^{\perp}_H \psi_1 \cdot \nabla_H \omega_1+\beta\psi_{1x}] \\
		&=-\nu\Delta \omega_2
		+q_2 [\nabla^{\perp}_H (2\psi_1-\psi_0 )\cdot \nabla_H (2\omega_1-\omega_0) +\beta (2\psi_1-\psi_0 )_x]  \\
		&+\nu\Delta (\omega_2-\omega_1)
		-(q_2-q_1)\nabla^{\perp}_H \psi_1 \cdot \nabla_H \omega_1 -q_2 \nabla^{\perp}_H \psi_1 \cdot \nabla_H (\omega_1-\omega_2)\\
		&-q_2\nabla^{\perp}_H (\psi_1-\psi_0) \cdot \nabla_H \omega_1
		-q_2 \nabla^{\perp}_H (\psi_1-\psi_0) \cdot \nabla_H (\omega_1-\omega_2)  \\
		&-\beta(q_2-q_1)\psi_{1x}-\beta q_2(\psi_1-\psi_0)_x.
		\end{split}
	\end{equation}
	Recall that  the invariant measure $\mu_k$ is supported on the $H^1$-uniformly bounded global attractor $\mathcal{A}_{q,k}$. By \eqref{cylin} and integration by parts, one has
	\begin{equation}\label{refo2}
		\begin{split}
		&\Big\langle \Delta (\omega_2-\omega_1), \Phi^\prime(\omega_1)\Big\rangle =\sum_{i=1}^m \frac{\partial \phi}{\partial z_i}  \Big\langle \Delta (\omega_2-\omega_1), g_i\Big\rangle \\
		&=\sum_{i=1}^m \frac{\partial \phi}{\partial z_i}  \Big\langle  (\omega_2-\omega_1), \lambda_i g_i\Big\rangle.
		\end{split}
	\end{equation}
    { In the following the integration is understood on the product space $(H^\alpha_p\times \mathbb{R})^2$ and with respect to $\left[ {V}^1,  {V}^0\right]^T$.}
	In light of the reformulation \eqref{refo1}, Eq. \eqref{refo2} and the asymptotic consistency \eqref{asym-con-eq},  we write the integrand in Eq. \eqref{Weak-form1} as
    \begin{align*}
        \Big\langle-f-\nu\Delta \omega_2 
		+q_2[ \nabla^{\perp}_H (2\psi_1-\psi_0 )\cdot \nabla_H (2\omega_1-\omega_0)+\beta(2\psi_1-\psi_0 )_x], \Phi^\prime\big(\frac{3}{2}\omega_1-\frac{1}{2}\omega_0\big)\Big\rangle,
    \end{align*}
which by the scheme translates into
\begin{align}
\label{IM-in2}
   -\Big\langle \frac{3}{2}\omega_2-\frac{1}{2}\omega_1-(\frac{3}{2}\omega_1-\frac{1}{2}\omega_0), \Phi^\prime\big(\frac{3}{2}\omega_1-\frac{1}{2}\omega_0\big)\Big\rangle. 
\end{align}
{ 
Define 
\[h(s):=\Phi\Big(\frac{3}{2}\omega_1-\frac{1}{2}\omega_0+s\big( \frac{3}{2}\omega_2-\frac{1}{2}\omega_1-\frac{3}{2}\omega_1+\frac{1}{2}\omega_0\big)\Big), \qquad s\in [0, 1].\]
By Taylor's theorem 
\[h(1)-h(0)=h^\prime(0)+\frac{1}{2}h^{\prime\prime}(\theta), \qquad \exists \theta \in [0, 1].\]
By the definition of cylindrical test functions and \eqref{cylin}, one has
\begin{align*}
&h^\prime(0)=\Big\langle \frac{3}{2}\omega_2-\frac{1}{2}\omega_1-(\frac{3}{2}\omega_1-\frac{1}{2}\omega_0), \Phi^\prime\big(\frac{3}{2}\omega_1-\frac{1}{2}\omega_0\big)\Big\rangle,\\
&h^{\prime\prime}(\theta)=\sum_{i=1}^m\sum_{j=1}^m \frac{\partial^2 \phi\big( \omega_\theta\big)}{\partial z_i \partial z_j} \Big(\frac{3}{2}\omega_2-\frac{1}{2}\omega_1-(\frac{3}{2}\omega_1-\frac{1}{2}\omega_0), g_i\Big) \times \Big(\frac{3}{2}\omega_2-\frac{1}{2}\omega_1-(\frac{3}{2}\omega_1-\frac{1}{2}\omega_0), g_j\Big),    
\end{align*}
where 
\[\omega_\theta:=\frac{3}{2}\omega_1-\frac{1}{2}\omega_0+\theta\big( \frac{3}{2}\omega_2-\frac{1}{2}\omega_1-\frac{3}{2}\omega_1+\frac{1}{2}\omega_0\big).\]
 By  the $H^{-1}$ asymptotic consistency estimate in \eqref{asym-con-eq}, the remainder is controlled by  
\begin{align*}
   &\leq C\sum_{i=1}^m\sum_{j=1}^m ||\frac{\partial^2 \phi}{\partial z_i \partial z_j}||_\infty \big(||\omega_2-\omega_1||_{H^{-1}}+||\omega_1-\omega_0||_{H^{-1}}\big)^2||\nabla g_i||\times ||\nabla g_j||\\
   &\leq C k^2\sum_{i=1}^m\sum_{j=1}^m ||\frac{\partial^2 \phi}{\partial z_i \partial z_j}||_\infty ||\nabla g_i||\times ||\nabla g_j||.
\end{align*}
    }
    Therefore \eqref{Weak-form1} becomes
	\begin{align*}
        &=-\lim_{k\rightarrow 0}\frac{1}{k}\int  \Phi\big(\frac{3}{2}\omega_2-\frac{1}{2}\omega_1\big)-\Phi\big(\frac{3}{2}\omega_1-\frac{1}{2}\omega_0\big)\, d \mu_k \\
		& =-\lim_{k\rightarrow 0}\frac{1}{k}\int  \Phi \left(\mathbb{P}_\omega\left\{\mathbb{S}_{q,k} \begin{bmatrix}
			{V}^1\\
			{V}^0
		\end{bmatrix} \cdot \begin{bmatrix}
			\frac{3}{2}\\
			-\frac{1}{2}
		\end{bmatrix}\right\}\right)- \Phi\left(\mathbb{P}_\omega\left\{\begin{bmatrix}
			{V}^1\\
			{V}^0
		\end{bmatrix} \cdot \begin{bmatrix}
			\frac{3}{2}\\
			-\frac{1}{2}
		\end{bmatrix}\right\}\right)\, d \mu_k \nonumber \\
		&=0,  
	\end{align*}
	where the last step follows from the  invariance of $\mu_k$ under $\mathbb{S}_{q,k}$. This establishes
	\begin{align}
		\label{Weak-form3}
		&\int_{H^\alpha_p} 
		\Big\langle-f-\nu \Delta \omega + \nabla^{\perp}_H \psi \cdot \nabla_H \omega+\beta\psi_x, \Phi^\prime(\omega)\Big\rangle\, d\mathbb{P}^\ast_\omega\mathbb{P}^\ast_1 \mu_0=0.
	\end{align}
	This is the weak form of invariance of $\mathbb{P}^\ast_\omega\mathbb{P}^\ast_1 \mu_0$ with respect to the  continuously stratified quasi-geostrophic equations \eqref{SVe}. 

       \noindent {  Step 3. Energy inequality.} 
       
	For the statistical version of the energy inequality for the PDEs \eqref{SVe}, one first establishes the corresponding one for the invariant measure $\mu_k$, and then takes the limit $k\rightarrow 0$. The $L^2$ energy estimate of the scheme \eqref{HDS-SVe1}--\eqref{HDS-SVe2} yields
	\begin{equation}
		\label{sta-en1}
		\begin{split}
		&\frac{1}{k}\left\{(||  V^{n+1}||_G^2+|Q^{n+1}|_G^2)- (||  V^{n}||_G^2+|Q^{n}|_G^2) \right\}
		+\nu  ||\nabla  \omega^{n+1}||^2\\
		&+\gamma |q^{n+1}|^2-(f, \omega^{n+1})-\gamma q^{n+1} \\
		&=-\frac{1}{k}\left\{\frac{1}{4}\big(||\omega^{n+1}-2\omega^n+\omega^{n-1}||^2+|q^{n+1}-2q^n+q^{n-1} |^2\big)\right\}.
		\end{split}
	\end{equation}
	Denoting $\widetilde{\mathbf{V}}=[{V}_1, {V}_0]^T$, and denoting by $(\cdot)_{1\omega}$ and $(\cdot)_{1q}$ the $\omega$- and $q$-components of the first coordinate, respectively, one derives from \eqref{sta-en1} and the discrete dynamical system formulation \eqref{dyna2} that
	\begin{align}
		\label{sta-en2}
		&0\geq \int   \frac{1}{k}\Big(||\mathbb{S}_{q,k} \widetilde{\mathbf{V}}||_G^2-|| \widetilde{\mathbf{V}}||_G^2\Big)+\nu  ||\nabla  \big(\mathbb{S}_{q,k} \widetilde{\mathbf{V}}\big)_{1\omega}||^2+\gamma |\big(\mathbb{S}_{q,k} \widetilde{\mathbf{V}}\big)_{1q}|^2 \nonumber \\
		&\quad \quad -\Big(f, \big(\mathbb{S}_{q,k} \widetilde{\mathbf{V}}\big)_{1\omega}\Big)-\gamma \big(\mathbb{S}_{q,k} \widetilde{\mathbf{V}}\big)_{1q}\, d\mu_k (\widetilde{\mathbf{V}}), \nonumber \\
		&=\int \nu  ||\nabla  \big(\mathbb{S}_{q,k} \widetilde{\mathbf{V}}\big)_{1\omega}||^2+\gamma |\big(\mathbb{S}_{q,k} \widetilde{\mathbf{V}}\big)_{1q}|^2-\Big(f, \big(\mathbb{S}_{q,k} \widetilde{\mathbf{V}}\big)_{1\omega}\Big)-\gamma \big(\mathbb{S}_{q,k} \widetilde{\mathbf{V}}\big)_{1q}\, d\mu_k (\widetilde{\mathbf{V}}), \nonumber \\
		&=\int \nu  ||\nabla  \omega_1||^2+\gamma | q_1 |^2-\big(f, \omega_1\big)-\gamma q_1\, d\mu_k (\widetilde{\mathbf{V}}),
	\end{align}
	where the invariance of $\mu_k$ under $\mathbb{S}_{q,k}$ is utilized in the last step. On the other hand, one recalls from \cite{Wang2012} that
	\begin{align*}
		\int \nu  ||\nabla  \omega_1||^2 \, d\mu_0 \leq  \liminf_{k\rightarrow 0} \int \nu  ||\nabla  \omega_1||^2 \, d\mu_k.
	\end{align*}
	Since $\mu_k$ is carried by the global attractor and on the global attractor $|q_1-1| \leq C k^{\frac{1}{2}}$,  taking $\liminf_{k\rightarrow 0}$ in \eqref{sta-en2} one obtains {  by the lower semi-continuity of the energy functional with respect to weak convergence of measures}
	\begin{align}
		\label{sta-en3}
		0&\geq \int_{(H^\alpha_p\times \mathbb{R})^2} \nu  ||\nabla  \omega_1||^2-\big(f, \omega_1\big)\, d\mu_0 (\widetilde{\mathbf{V}}), \nonumber \\
		&=\int_{H^\alpha} \nu  ||\nabla  \omega||^2-\big(f, \omega\big)\, d\mathbb{P}^\ast_\omega\mathbb{P}^\ast_1 \mu_0(\widetilde{\mathbf{V}}).
	\end{align}
	In light of the identity \eqref{Weak-form3} and the energy inequality \eqref{sta-en3}, one concludes that $\mathbb{P}^\ast_\omega\mathbb{P}_1^\ast \mu_0 \in \mathcal{IM}$. This completes the proof.
\end{proof}

\section{Numerical experiments}\label{s-nu}

{ In this section we present numerical experiments illustrating the main properties of the proposed mr-SAV-BDF2 scheme. The purpose of the experiments is threefold. First, manufactured-solution tests confirm the expected second-order temporal accuracy. Second, Kolmogorov-flow simulations demonstrate the long-time robustness of the method for a forced dissipative flow with nontrivial dynamics. Third, we examine intermittent bursting statistics through cumulative probabilities of selected observables under time-step refinement. These numerical tests are intended to complement, rather than replace, the theoretical convergence results for attractors and invariant measures.}

Since the 2D NSE in vorticity stream function formulation is well studied, most numerical experiments are performed solving the 2D NSE under periodic boundary conditions, with only a brief mention of the result for the continuously stratified QG model. Since periodic boundary conditions are used, the mean-reverting SAV-BDF2 scheme in the stream function formulation is  implemented in Matlab utilizing Fourier collocation methods.
In all tests, the scheme is initiated by a first-order version of the mr-SAV-BDF2 scheme.

\subsection{ Accuracy} 
We perform tests for verification of accuracy using manufactured solutions.  
For the 2D NSE in the vorticity stream function formulation, the exact solutions are 
\begin{align*}
	&  \psi=\cos(t) \cos(2\pi x) \cos(2\pi y), \quad \omega=-\Delta \psi.
\end{align*}
The parameters are $Re=10, \gamma=1000, T=100, \Omega=(0, 1)^2$. $32$ Fourier modes are employed. The relative error in the $l^\infty$ norm and convergence order is displayed in Table \ref{stream_conv}, in which the short-hand notation pre-asymp means pre-asymptotic convergence order before the mesh becomes sufficiently fine.
\begin{table}[hbt]
	\begin{center}
		\caption{Relative error   for    $\omega$ and $\psi$  in the $l^\infty$ norm for the NSE in the vorticity stream function formulation,   with  $32$ Fourier modes in each direction} \label{stream_conv}
		\begin{tabular}{c c c c c}
			\hline  $k$  &  $\|e_{\omega}\|_{\infty}$    &order  &      $\|e_\psi\|_{\infty}$    &order    \\[0.5ex]\hline
			0.05    &1.36619   &       &1.818292e-03   &                                \\
			0.025    & 1.025531e-03  & pre-asymp        &8.159268e-06   & pre-asymp   \\
			0.0125   & 2.553669e-04  &2.00      &2.031416e-06  &2.00  \\
			0.00625  & 6.363869e-05  &2.00       &5.064025e-07   &2.00   \\
			0.003125   & 1.588605e-05  &2.00     &1.264143e-07  &2.00  \\
			0.0015625  & 3.968483e-06  &   2.00        & 3.157964e-08  &  2.00                                    \\
			0.00078125 & 9.917372e-07 &    2.00       &7.891879e-09    & 2.00     
			\\			\hline
		\end{tabular}
	\end{center}
\end{table}	

We take similar manufactured solution for the continuously stratified QG model, 
The domain is a unit cube $[0, 1]^3$, $Re=10$, $\beta=F=1$, $T=100$, $\gamma=1000$. Again 32 Fourier modes are utilized. Table \ref{conv_CQG} displays the relative error in the $l^\infty$ norm and convergence order.
\begin{table}[hbt]
	\begin{center}
		\caption{Relative error   for    $\omega$ and $\psi$  in the $l^\infty$ norm for the continuously stratified QG model,   with  $32$ Fourier modes in each direction} \label{conv_CQG}
		\begin{tabular}{c c c c c}
			\hline  $k$  &  $\|e_{\omega}\|_{\infty}$    &order  &      $\|e_\psi\|_{\infty}$    &order    \\[0.5ex]\hline
			0.05    &5.040908e-05   &       &5.040795e-05   &                                \\
			0.025    & 1.231334e-05  &2.03        &1.231334e-05   &2.03  \\
			0.0125   & 3.041844e-06  &2.01      &3.041844e-06  &2.01  \\
			0.00625  & 7.558699e-07  &2.00       &7.558699e-07   &2.00   \\
			0.003125   & 1.883918e-07  &2.00     &1.883918e-07  &2.00  \\
			0.0015625  & 4.702588e-08  &   2.00        & 4.702588e-08  &  2.00                                    \\
			0.00078125 & 1.174745e-08 &    2.00       &1.174745e-08    & 2.00     
			\\			\hline
		\end{tabular}
	\end{center}
\end{table}

\subsection{Robustness and long-time stability}	
As a benchmark test, one considers the 2D NSE with the Kolmogorov forcing and periodic boundary conditions in the domain $\Omega=(0, 2\pi)^2$, as is detailed in \cite{Armbruster1996}. A typical Kolmogorov forcing takes the form
\begin{align}\label{Kol}
	\mathbf{f}:=\Big[\frac{m^3}{Re} \cos(my), 0\Big]^T, \quad m \in \mathbb{N}^+.
\end{align}
The steady-state solution, also known as the basic Kolmogorov flow, is given in a streamfunction formulation as $\psi=\sin(my)$.  The steady-state solution is stable when $Re$ is relatively small. It becomes unstable as $Re$ crosses a critical value and coherent structure develops. 

{ 
\subsubsection{Robustness} 
We examine the robustness of the mr-SAV-BDF2 scheme in terms of the closeness between $q_n$ and the desired value $1$ in the long-time simulation of the 2D NSE in the form of the vorticity-stream function. The parameters are $Re=100$, $k=0.002$, $T=10000$, and the force is of the Kolmogorov type. Fig. \ref{mr_plot100} compares the first 100 unit time of $|q_n-1|$ generated by the SAV-BDF2 scheme without mean-reverting $\gamma=0$ (left panel), and by the mr-SAV-BDF2 scheme with $\gamma=1000$ (right panel). The error of the non-mean-reverting SAV method increases to the order of $10^{-3}$ and stays elevated, while that of the mr-SAV-BDF2 scheme peaks at about $10^{-6}$ and then generally decreases though oscillating. 
Fig. \ref{mr_plot10000} compares the performance in the asymptotic regime $t\in[100, 10000]$. It shows that the error grows linearly in time, though slowly, in the non-mean-reverting version of the SAV-BDF2 method (left panel), whereas that of the mr-SAV-BDF2 scheme oscillates between 0 and $10^{-8}$ (right panel). This numerical experiment demonstrates that the mr-SAV-BDF2 method is more robust  than the non-mean-reverting SAV method in the long-time computation.
\begin{figure}[!ht]
	\subfigure
	{
		\centering
		\begin{minipage}[t]{1.0\linewidth}
			\includegraphics[width=2.6in]{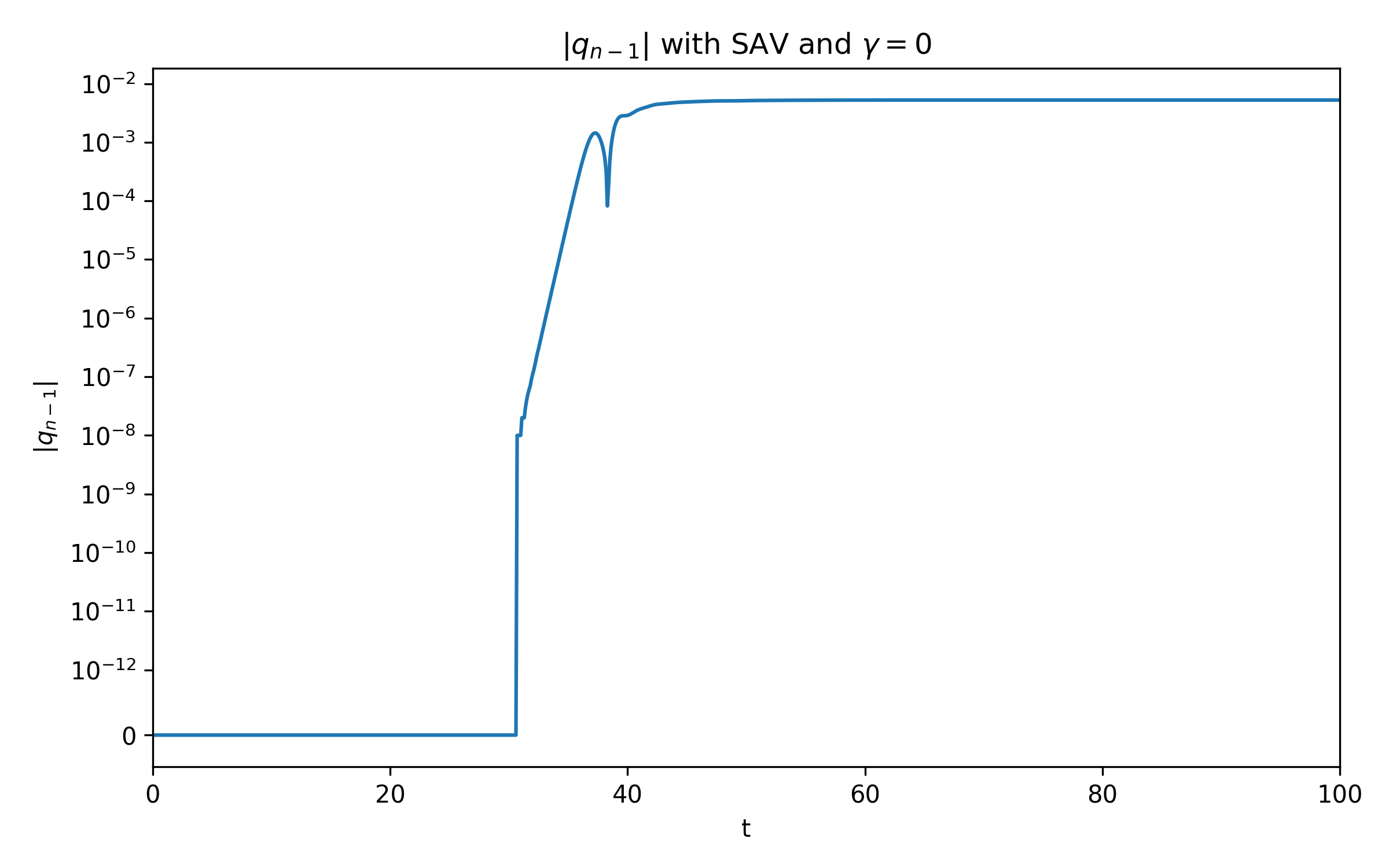}
			\includegraphics[width=2.6in]{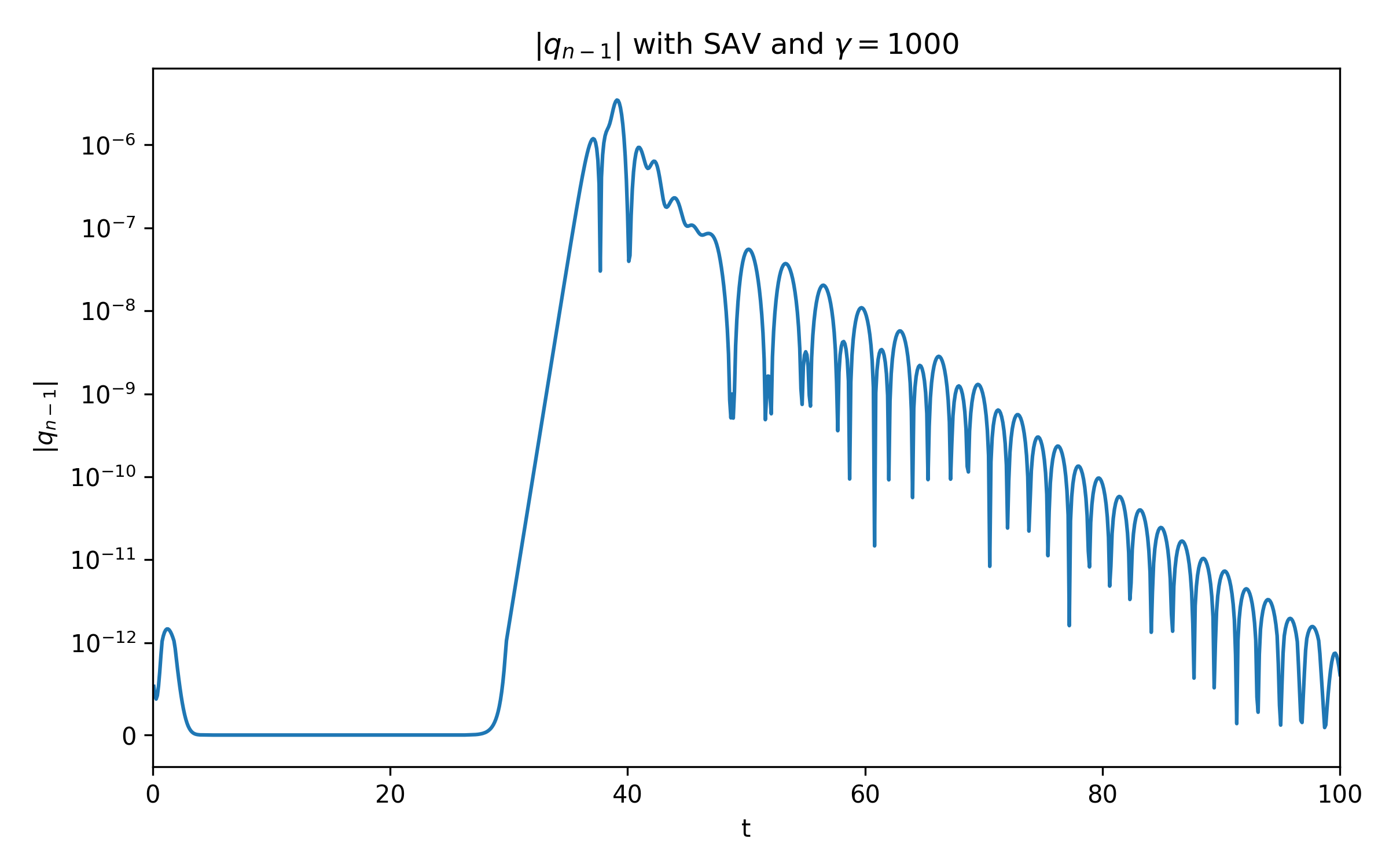}
		\end{minipage}
	}
	\caption{The evolution of $|q_n-1|$ for $t\in [0, 100]$ by the SAV-BDF2 scheme without mean-reverting $\gamma=0$ (left),   and  with mean-reverting $\gamma=1000$ (right), respectively. }
	\label{mr_plot100}
\end{figure}

\begin{figure}[!ht]
	\subfigure
	{
		\centering
		\begin{minipage}[t]{1.0\linewidth}
			\includegraphics[width=2.6in]{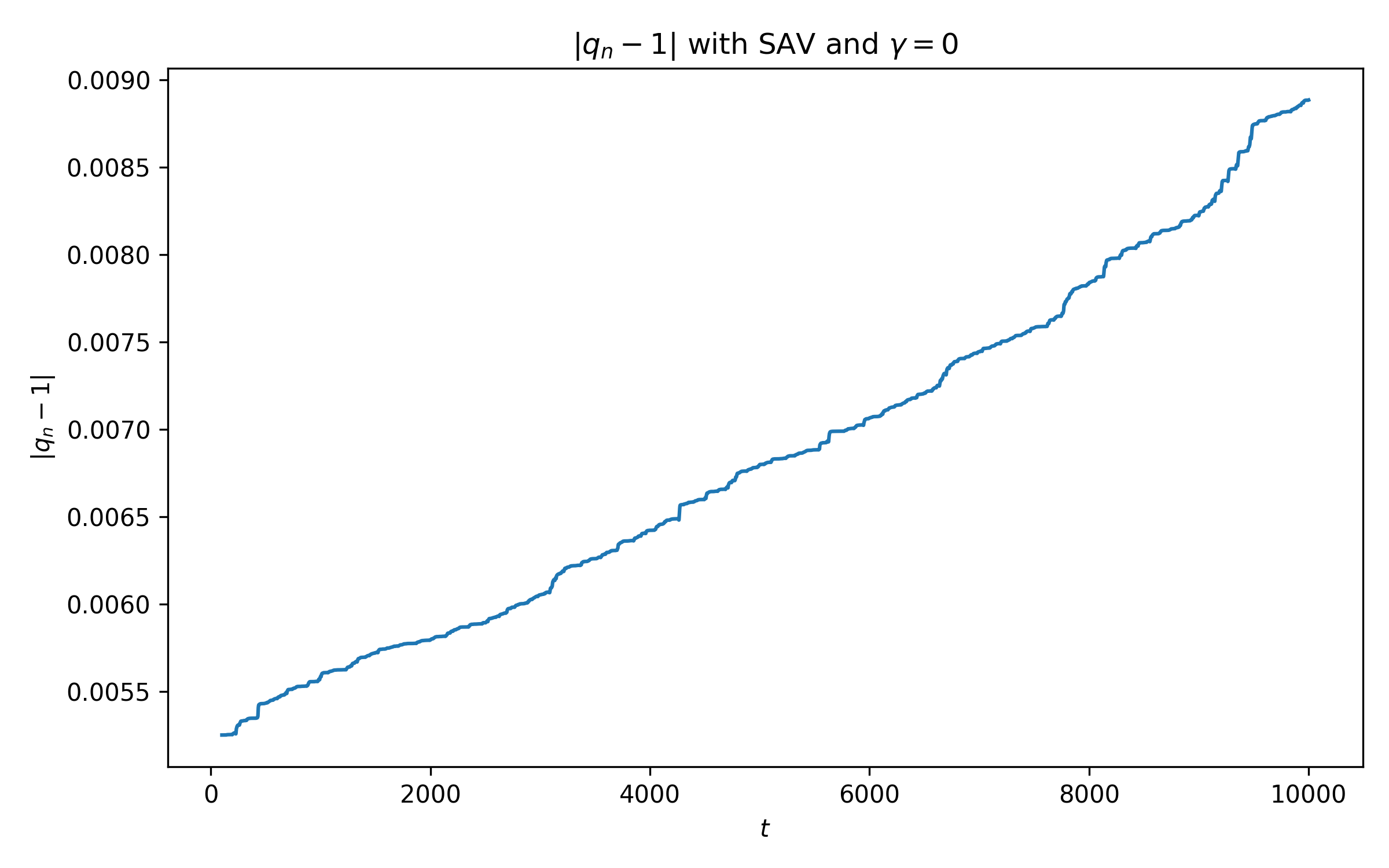}
			\includegraphics[width=2.6in]{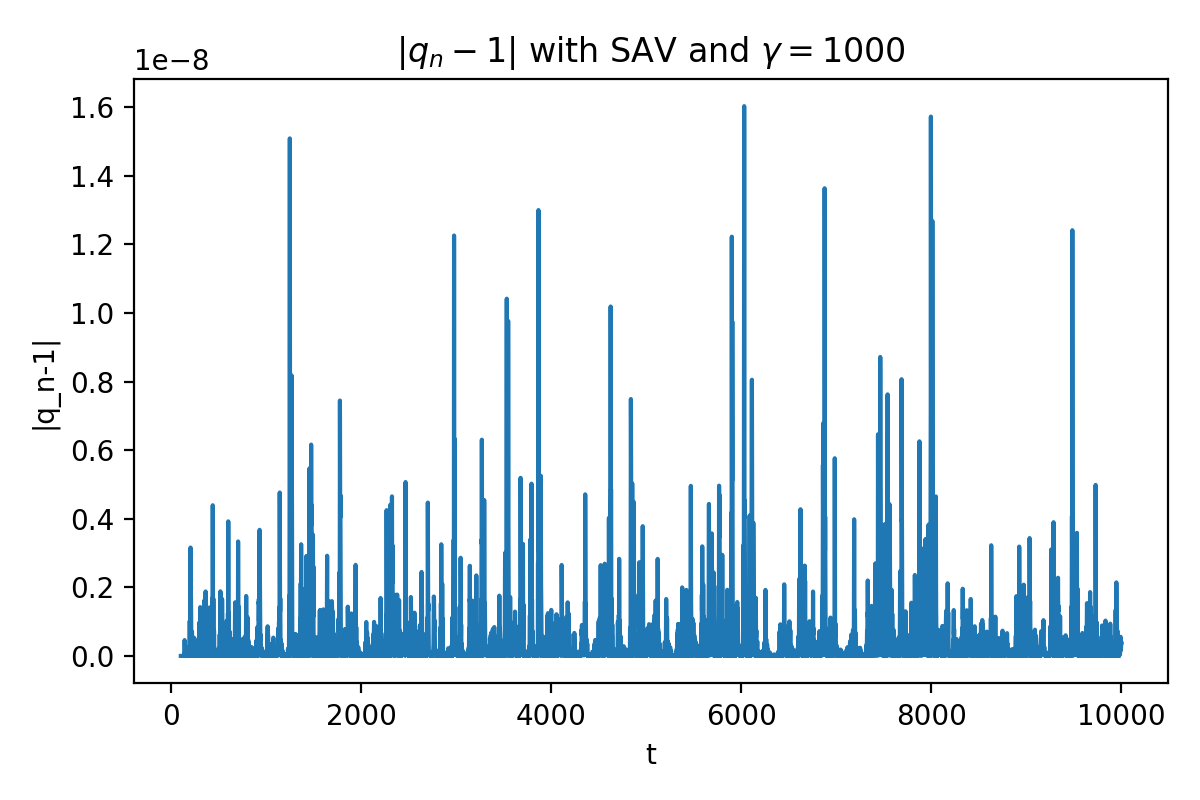}
		\end{minipage}
	}
	\caption{The evolution of $|q_n-1|$ for $t\in [100, 10000]$ by the SAV-BDF2 scheme without mean-reverting $\gamma=0$ (left),   and  with mean-reverting $\gamma=1000$ (right), respectively. }
	\label{mr_plot10000}
\end{figure}

}

\subsubsection{Long-time stability}
One verifies the long-time stability of the mr-SAV-BDF2 scheme for the 2D NSE  in the stream function formulation with the  Kolmogorov forcing \eqref{Kol} and $m=2, Re=100$. The initial stream function is a (periodic) perturbation of the steady-state solution
\begin{align*}
	\psi(0, x, y)=\sin(2y)+0.001\sin(2 x)\sin(2 y).
\end{align*}
The evolution of the $L^2$ norm and $H^1$ norm of the vorticity by the mr-SAV-BDF2 scheme for $T=1000$ is depicted in Fig. \ref{Sta_plot}  with 256 Fourier modes and $k=0.01, 0.005, 0.0025$ respectively. It is noted that the BDF2 IMEX scheme with explicit discretization of the advection term blows up for $k=0.003$, consistent with the time-step requirement for long-time stability of this BDF2-Gear's extrapolation scheme \cite{Wang2012, ChWa2016}.

\begin{figure}[!ht]
	\subfigure
	{
		\centering
		\begin{minipage}[t]{1.0\linewidth}
			\includegraphics[width=2.5in]{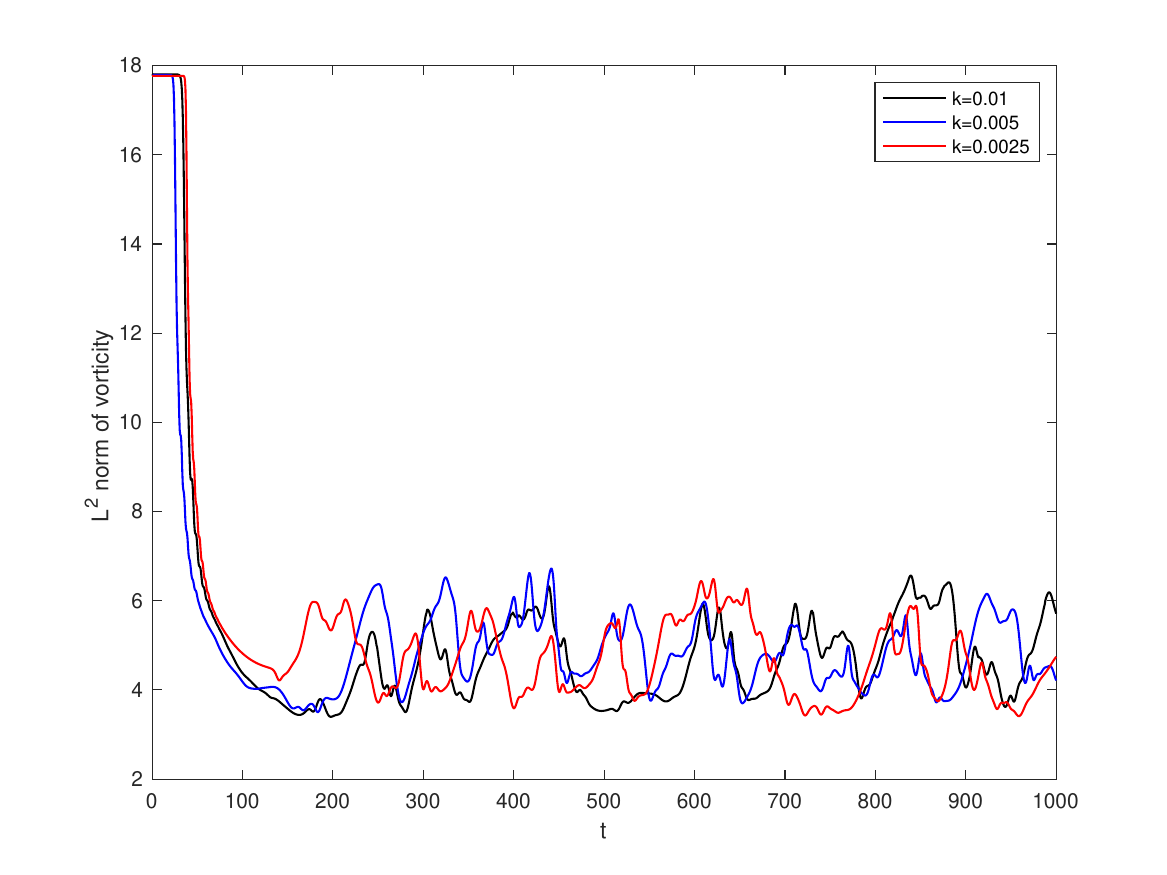}
			\includegraphics[width=2.5in]{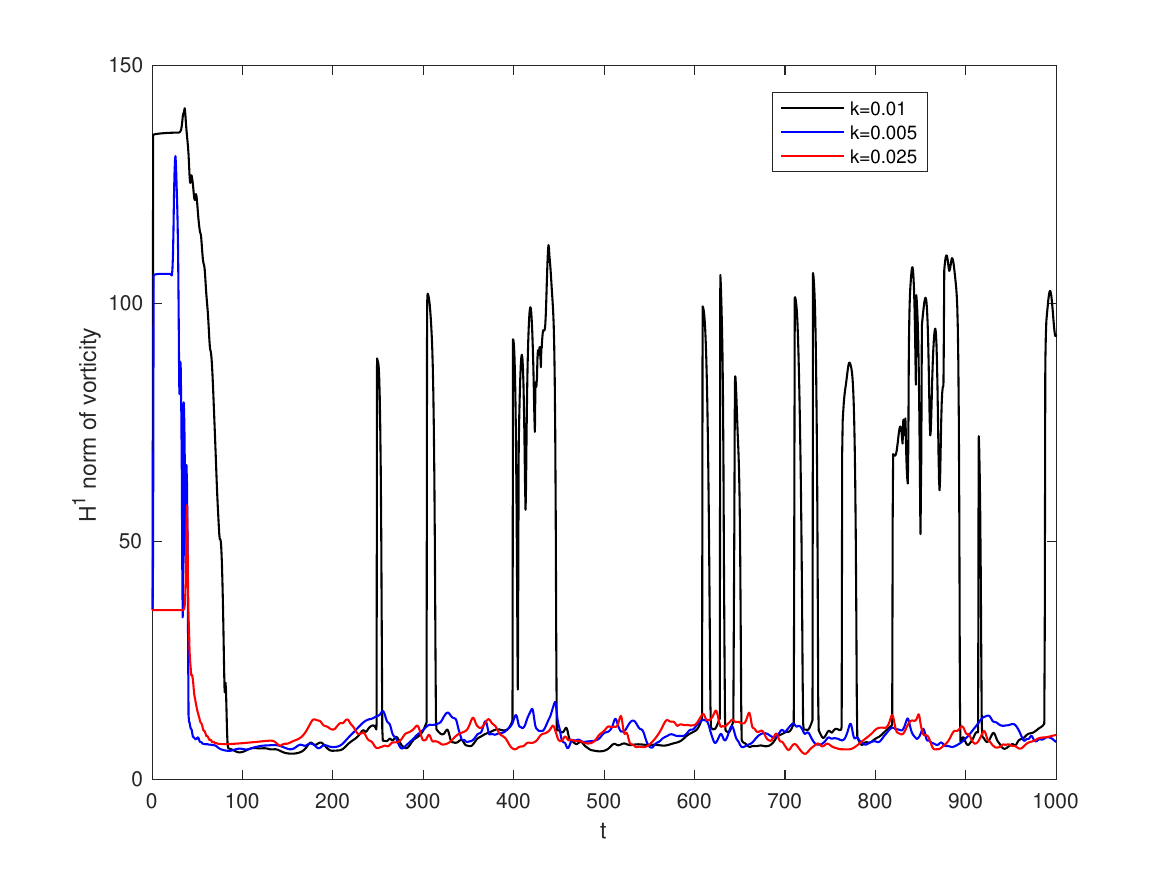}
		\end{minipage}
	}
	\caption{The $L^2$ norm (enstrophy) and $H^1$ norm (palinstrophy) of the vorticity as a function of time by the mr-SAV-BDF2 scheme with 256 Fourier modes and $k=0.01, 0.005, 0.0025$ respectively.}
	\label{Sta_plot}
\end{figure}

\subsection{Statistical persistence of intermittent events}
{ Since the exact invariant measure is unavailable, we do not attempt to measure the distance between invariant measures directly. 
Instead, we compute the empirical cumulative probabilities associated with bursting events. 
The stability of these probabilities as $k$ is refined provides an observable-level consistency check for the long-time statistics.}

The initial condition used in this simulation is a perturbation of the basic Kolmogorov flow
\begin{align*}
	\psi(0, x, y)=\sin(2y)+0.001\sin(2\pi x)\sin(2\pi y).
\end{align*}
It is discovered through numerical simulation in \cite{Armbruster1996} that at $Re=25.70$ the solution is quasi-periodic, consisting of a traveling structure plus additional time-dependent behavior, and  bursts occur intermittently at $Re \geq 25.77$. In the following simulation, the Reynolds number is taken to be $25.7715$, $k=0.001$, and the final time is $T=10000$. Fourier collocation method is utilized for spatial discretization with $256$ Fourier modes. {  The eddy turnover time in this test is about $T_e=\pi/2\approx1.6$. Treating the first $1000$ unit computation time as transient, the total usable statistical sampling time is $9000/1.6 \approx 5625 T_e$. Hence this is roughly sufficient for a preliminary persistence test of exceedance probability. }

{ The real part of the Fourier coefficient of mode $e^{iy}$ and the palinstrophy, i.e., the $L^2$ norm of the gradient, are plotted in Fig. \ref{mode-L2} as a function of time.} Bursts occur intermittently characterized by a sudden dramatic change of magnitude of the Fourier coefficient.  The dynamics undergo a long laminar regime, then a chaotic explosion (burst) ensues, followed by another long period of laminar regime. The pattern repeats in time. {  In particular, we note that the computation time has to be long enough ($T>1000$  in this case) in order to capture the intermittent bursting events, underscoring the necessity of long-time stability (uniform-in-time estimate) and efficiency of numerical methods.}


\begin{figure}[!ht]
	\subfigure
	{
		\centering
		\begin{minipage}[t]{1.0\linewidth}
			\includegraphics[width=2.6in]{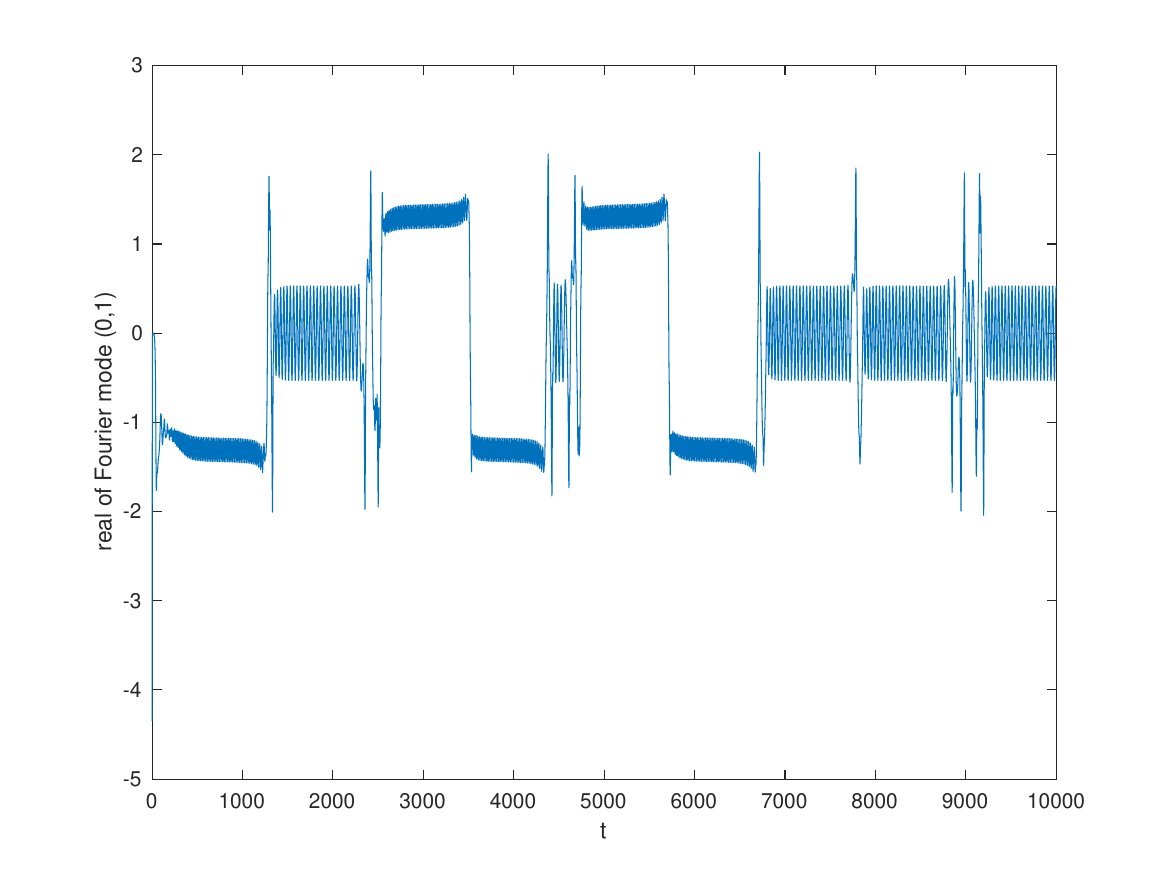}
			\includegraphics[width=2.6in]{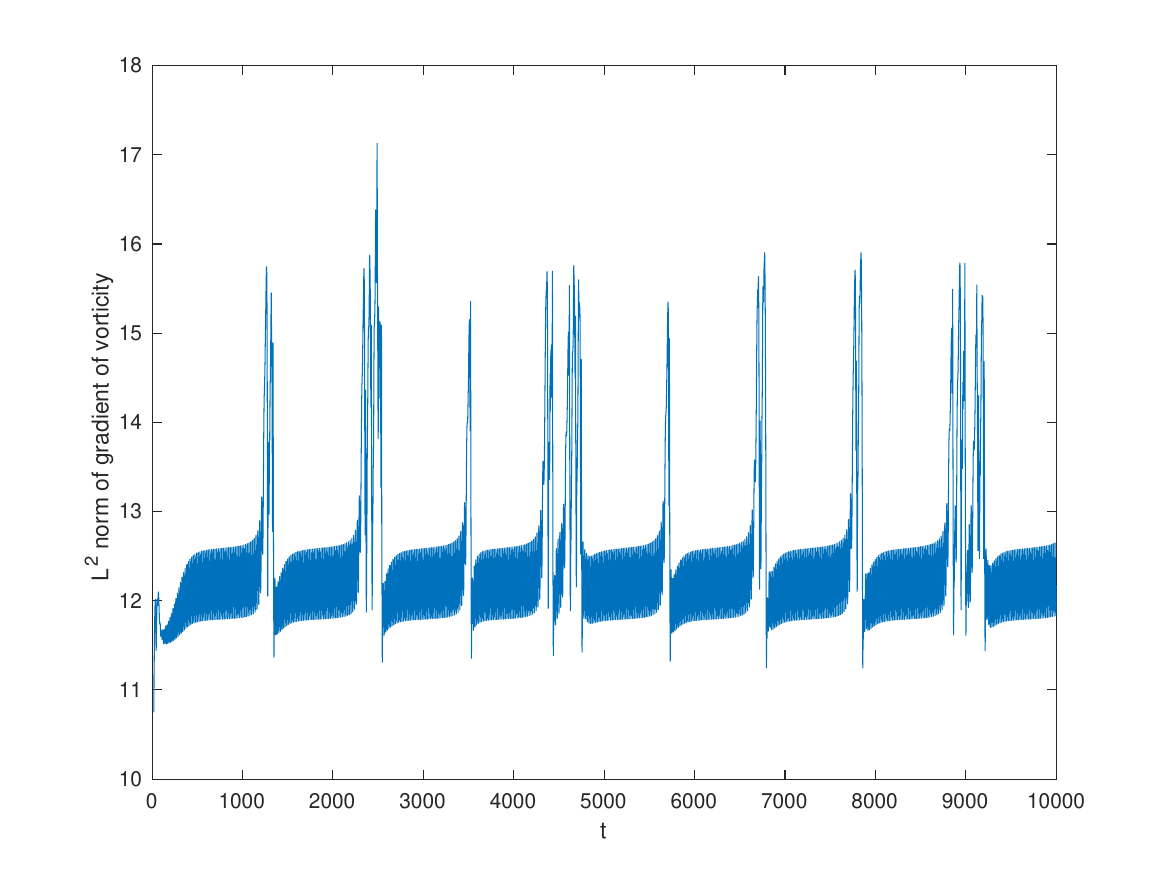}
		\end{minipage}
	}
	\caption{Left panel: the real part of the Fourier coefficient of mode $e^{iy}$ as a function of time; right panel: the palinstrophy ($L^2$ norm of  gradient of vorticity) as a function of time. $Re=25.7715, k=0.001$ with 256 Fourier modes.}
	\label{mode-L2}
\end{figure}

{ 
Table \ref{tab:max_vorticity_statistics} reports the empirical mean, standard deviation, and normalized tail probabilities of the maximum vorticity over long-time simulations with several time-step sizes. The results should be interpreted qualitatively: the purpose is to verify that the intermittent bursting regime persists under refinement, rather than to provide high-precision estimates of rare-event probabilities. The far-tail probabilities are sensitive to the effective number of statistically independent bursts, and a dedicated rare-event study would require longer trajectories or block-based uncertainty estimates.
\begin{table}[ht]
\centering
\caption{Post-transient statistics of the maximum vorticity and exceedance
probabilities for two burst thresholds.}
\begin{tabular}{ccccc}
\hline
$\Delta t$
& Mean $\mu$
& Std.\ Dev.\ $\sigma$
& $P(X>\mu+2\sigma)$
& $P(X>\mu+2.5\sigma)$ \\
\hline
0.005     & 3.141004 & 0.393018 & 0.068980 & 0.046083 \\
0.0025    & 3.106066 & 0.399789 & 0.063278 & 0.045549 \\
0.00125   & 3.186472 & 0.450708 & 0.076718 & 0.046489 \\
0.000625  & 3.220583 & 0.485382 & 0.065535 & 0.021166 \\
0.0003125 & 3.120307 & 0.408605 & 0.063044 & 0.045155 \\
\hline
\end{tabular}
\label{tab:max_vorticity_statistics}
\end{table}
}

\section{Conclusion} \label{conc}
A novel mean-reverting scalar auxiliary variable BDF2 (mr-SAV-BDF2) method  is introduced to preserve the underlying dissipative structure of the continuous dynamical system that yields a uniform-in-time estimate of the numerical solution. For geophysical fluid models such as the damped-driven barotropic quasi-geostrophic model, the 2D Navier-Stokes equations in vorticity--stream function formulation, as well as the damped-driven continuously stratified QG model, we are also able to derive a uniform large-time estimate of the vorticity in $H^1$ norm in addition to the $L^2$ norm guaranteed by the general framework. In the case of  time-independent forcing, the numerical scheme defines an autonomous discrete dissipative dynamical system for a suitable choice of phase space. We show that the global attractor and the invariant measures (stationary statistical solutions) of the discrete dynamical system converge to those of the continuous models at vanishing time-step size.
Numerical results { are provided to gauge} the performance of the new algorithm in terms of accuracy, efficiency, stability and robustness.  For the 2D NSE with the no-slip no-penetration boundary condition, the unconditional uniform $H^1$ estimate remains open.

\section*{Appendix  I}
{ 
Here we provide the proof of the technical  Lemma \ref{Lem-key}.
\begin{proof}[Proof of Lemma~\ref{Lem-key}]
Set
\[
\rho_k:=1+bl_k,
\]
and define the augmented energy
\begin{equation}
E_n
:=
Y_n+\mu k\left(\frac14+\varepsilon\right)r_n
+\varepsilon\mu k r_{n-1}.
\label{Aug_E}
\end{equation}
We write the inequality \eqref{Bas_ine} as
\[
Y_{n+1}+\mu k r_{n+1}
\le
Y_n+\varepsilon\mu k r_n+\varepsilon\mu k r_{n-1}+kg.
\]
Adding $\frac{\mu}{4}kr_n$ to both sides gives
\begin{equation}
Y_{n+1}+\mu k r_{n+1}+\frac{\mu}{4}kr_n
\le
E_n+kg.
\label{2}
\end{equation}

We next show that the left-hand side of \eqref{2} controls $\rho_kE_{n+1}$.
Write
\[
\mu k r_{n+1}+\frac{\mu}{4}kr_n
=
\frac{\mu}{24}k(r_{n+1}+r_n)
+
\frac{\mu}{24}k(r_{n+1}+r_n)
+
\frac{11\mu}{12}kr_{n+1}
+
\frac{\mu}{6}kr_n .
\]
Hence
\begin{align}
&Y_{n+1}+\mu k r_{n+1}+\frac{\mu}{4}kr_n
\nonumber\\
&=
\left[
Y_{n+1}
+
\frac{\mu}{24}k(r_{n+1}+r_n)
\right]
+
\frac{11\mu}{12}kr_{n+1}
+
\frac{\mu}{6}kr_n
+
\frac{\mu}{24}k(r_{n+1}+r_n).
\label{3}
\end{align}
Since $k\ge l_k$, $b\le \mu c_d/24$ and in light of \eqref{Bas_as}, we obtain
\[
\frac{\mu}{24}k(r_{n+1}+r_n)
\ge
\frac{\mu c_d}{24}kY_{n+1}
\ge
bl_kY_{n+1}.
\]
Therefore
\begin{equation}
Y_{n+1}
+
\frac{\mu}{24}k(r_{n+1}+r_n)
\ge
(1+bl_k)Y_{n+1}
=
\rho_kY_{n+1}.
\label{4}
\end{equation}
Moreover, since $b\le 1$ and $l_k\le 1$, we have $\rho_k\le 2$.
Using $0\le \varepsilon\le 1/12$,
\[
\rho_k\mu\left(\frac14+\varepsilon\right)
\le
2\mu\left(\frac14+\frac1{12}\right)
=
\frac{2\mu}{3}
<
\frac{11\mu}{12}, \qquad \rho_k\varepsilon\mu
\le
\frac{\mu}{6}.
\]
Thus
\begin{equation}
\frac{11\mu}{12}kr_{n+1}
+
\frac{\mu}{6}kr_n
\ge
\rho_k\mu k\left(\frac14+\varepsilon\right)r_{n+1}
+
\rho_k\varepsilon\mu k r_n .
\label{5}
\end{equation}
Combining \eqref{3}, \eqref{4} and \eqref{5}, we get
\[
Y_{n+1}+\mu k r_{n+1}+\frac{\mu}{4}kr_n
\ge
\rho_kE_{n+1}
+
\frac{\mu}{24}k(r_{n+1}+r_n).
\]
Together with \eqref{2}, this yields
\begin{equation}
\rho_kE_{n+1}
+
\frac{\mu}{24}k(r_{n+1}+r_n)
\le
E_n+kg.
\label{6}
\end{equation}

Dividing \eqref{6} by $\rho_k$ and iterating from $M$ to $n$ gives
\begin{equation}
E_{n+1}
+
\frac{\mu}{24}k
\sum_{i=M}^{n}
\frac{r_{i+1}+r_i}{\rho_k^{\,n+1-i}}
\le
\frac{E_M}{\rho_k^{\,n+1-M}}
+
kg
\sum_{i=M}^{n}
\frac{1}{\rho_k^{\,n+1-i}}\leq 
\frac{E_M}{\rho_k^{\,n+1-M}}+\frac{kg}{bl_k}.
\label{7}
\end{equation}

In view of the definition of $l_k$ and $E_M$,  \eqref{7} immediately yields for $k\leq 1$,
\begin{align}
    \label{kl1}
    Y_{n+1} \leq \frac{
Y_M+\mu \left(\frac14+\varepsilon\right)r_M
+\varepsilon\mu  r_{M-1}
}{
(1+bl_k)^{\,n+1-M}
}+\frac{g}{b}.
\end{align}
On the other hand, for $k>1$ keeping only the last term in the sum in \eqref{7}, and using \eqref{Bas_as} we obtain
\begin{align}
    \label{kg1}
    Y_{n+1} \leq \frac{48}{\mu c_d}\frac{
Y_M+\mu \left(\frac14+\varepsilon\right)r_M
+\varepsilon\mu  r_{M-1}
}{
(1+bl_k)^{\,n+1-M}
}+\frac{48}{\mu c_d}\frac{g}{b}.
\end{align}

Combining \eqref{kl1} and \eqref{kg1} gives the desired estimate \eqref{lem_fr}. This completes the proof. 
\end{proof}
}

\section*{Appendix II}
The purpose of this appendix is to demonstrate the convergence of the numerical scheme to that of the extended system for initial data from the global attractor of the scheme. This is needed in the convergence of the global attractors.
This convergence analysis differs from the standard ones in at least two aspects: (1) the convergence is uniform with respect to initial data from the global attractors of the scheme; (2) the initial data may not satisfy the second-order consistency requirement, i.e., $\omega^{1}-\omega^0$ is not of the order of $k^2$.

\begin{proof}[Proof of Lemma \ref{Traj-convL}]
\noindent {Step 1.}
{ We first recall the following regularity result for the damped-driven continuously stratified QG model.}  
For any $ \left[(\omega^1, q^1), (\omega^0, q^0)\right]^T \in \mathcal{A}_{q, k}$ and $f\in L^2(\Omega)$, one derives from Theorem \ref{main-co} and elementary energy estimates that  the solution to the expanded PDE system  \eqref{SVe1}--\eqref{SVe2} with initial data $(\omega^0, q^0) \in H^1_{p0}\times \mathbb{R}$ has the following regularity
	\begin{align}
		\label{Reg-Tra}
		\sup_{0\leq t \leq T} \big[||\nabla \omega||^2+||\omega_t||_{H^{-1}}^2\big]+\int_0^T ||\omega_t||^2 + ||\omega_{tt}||_{H^{-2}}^2\, dt \leq K.
	\end{align}

	
	The invariance of the global attractor $\mathcal{A}_{q,k}$ under $\mathbb{S}_{q, k}$ and the uniform boundedness \eqref{att-boun} implies
	\begin{align}
		||\nabla \omega^n|| \leq C(R_0, R_1), \quad \forall n\geq0.
	\end{align}

\noindent{Step 2.} { We are now ready to write down the error equation.}
	Adopting the notation
	\begin{align*}
		e_\omega^n:=\omega^n-\omega(t^n, \cdot), \quad \eVO^{n}=\big[e_\omega^n, e_\omega^{n-1}\big]^T,
	\end{align*}
	one derives the error equation
	\begin{align}
		\delta_t \eo^{n+1} -\nu \Delta  \eo^{n+1}&=\nabla^{\perp}_H\psi\big(t^{n+1}\big)\cdot \nabla_H\omega\big(t^{n+1}\big) -q^{n+1} \nabla^{\perp}_H \overline{\psi}^{n+1} \cdot \nabla_H \overline{\omega}^{n+1}\nonumber \\
		&\quad +\beta\big[\psi_x\big(t^{n+1}\big)-q^{n+1}\overline{\psi}^{n+1}_x\big]+R_{\omega}^{n+1}, \label{Err-SV1}
	\end{align}    
	where the consistency error function is defined as
	\begin{equation}
		\label{err-con}
		\begin{split}
		R_{\omega}^{n+1}&=\frac{3\omega \big(t^{n+1}\big)-4\omega \big(t^{n}\big)+\omega \big(t^{n-1}\big)}{2k}-\omega_t\big(t^{n+1}\big)  \\
		&=\frac{\omega \big(t^{n+1}\big)-\omega \big(t^{n}\big)}{k}-\omega_t\big(t^{n+1}\big)+\frac{\omega \big(t^{n+1}\big)-2\omega \big(t^{n}\big)+\omega \big(t^{n-1}\big)}{2k} \\
		&=\frac{1}{k}\int_{t^n}^{t^{n+1}} \omega_{tt}(s)(s-t^n)\, ds\\
		&+\frac{1}{2k}\Big[ \int_{t^n}^{t^{n+1}} \omega_{tt}(s)(t^{n+1}-s)\, ds+\int_{t^{n-1}}^{t^{n}} \omega_{tt}(s)(s-t^{n-1})\, ds  \Big].
		\end{split}
	\end{equation}
{  Note that the consistency error term is expressed in terms of the second-order time derivative, in light of the regularity estimate \eqref{Reg-Tra}.}
	By Minkowski's inequality and the Cauchy-Schwarz inequality, one derives from \eqref{err-con} 
	\begin{align}
		&||\Delta^{-1} R_\omega^{n+1}||^2 \leq Ck \int_{t^{n-1}}^{t^{n+1}}||\Delta^{-1}\omega_{tt}||^2\, ds. \label{con-es1}
	\end{align}

\noindent {Step 3.} {  The error estimate is then somewhat standard provided we track the dependence on the norm of the initial data and the final time $T$.
}

	Denoting $\nabla^{-1}=(-\Delta)^{-\frac{1}{2}}$ and testing Eq. \eqref{Err-SV1} with $-k\Delta^{-1}\eo^{n+1}$ and performing integration by parts, one obtains
	\begin{equation}
		\label{Err-SV-es-1}
		\begin{split}
		&||\nabla^{-1}\eVO^{n+1}||_{G}^2-||\nabla^{-1}\eVO^{n}||_{G}^2+k \nu ||e_\omega^{n+1}||^2 \\
		&\leq k\Big|\Big(\nabla^{\perp}_H\big[\psi\big(t^{n+1}\big)-\overline{\psi}^{n+1}\big]\cdot \nabla_H \Delta^{-1}\eo^{n+1},  \omega\big(t^{n+1}\big)\Big)\Big| \\
		&+k\Big|\Big(\nabla^{\perp}_H\overline{\psi}^{n+1}\cdot \nabla_H \Delta^{-1}\eo^{n+1},  \big[\omega\big(t^{n+1}\big)-\overline{\omega}^{n+1}\big] \Big)\Big| \\
		& +k \big|q^{n+1}-1\big|\left\{\Big|\Big(\nabla^{\perp}_H \overline{\psi}^{n+1} \cdot \nabla_H \Delta^{-1}\eo^{n+1},  \overline{\omega}^{n+1}\Big)\Big| 
		+ \beta\Big|\Big(\overline{\psi}_x^{n+1}, \Delta^{-1}\eo^{n+1}  \Big)\Big|\right\} \\
		& +k \beta\Big|\Big(\psi_x\big(t^{n+1}\big)-\overline{\psi}_x^{n+1},\Delta^{-1}\eo^{n+1}\Big)\Big|+k \Big|\Big(R_{\omega}^{n+1},   \Delta^{-1}\eo^{n+1}\Big)\Big|\\
		&:= \sum_{i=1}^5 I_i. 
		\end{split}
	\end{equation}
	Define
	\begin{align}
		&D \psi\big(t^{n+1}\big)=\psi\big(t^{n+1}\big)-\psi\big(t^{n}\big)-\big[\psi\big(t^{n}\big)-\psi\big(t^{n-1}\big)\big], \label{diff-con1}\\
		&D \omega\big(t^{n+1}\big)=\omega\big(t^{n+1}\big)-\omega\big(t^{n}\big)-\big[\omega\big(t^{n}\big)-\omega\big(t^{n-1}\big)\big] \label{diff-con2}\\
		&=\int_{t^n}^{t^{n+1}}\omega_t(s)\, ds-\int_{t^{n-1}}^{t^{n}}\omega_t(s)\, ds. \nonumber 
	\end{align}
    	Then
	\begin{align*}
		&\psi\big(t^{n+1}\big)-\overline{\psi}^{n+1}=D \psi\big(t^{n+1}\big)-2e_\psi^n+e_{\psi}^{n-1}.   
	\end{align*}
    	It follows from Minkowski's inequality and the Cauchy-Schwarz inequality that
	\begin{align}
		&||D \omega \big(t^{n+1}\big)||^2 \leq Ck \int_{t^{n-1}}^{t^{n+1}}||\omega_{t}||^2\, ds. \label{con-es2}
	\end{align}
	One derives
	\begin{align}
		\label{non-Tra-1}
		I_1 &\leq C(\epsilon, \nu)k||\nabla \omega\big(t^{n+1}\big)||^2 ||\nabla^{\perp}_H\big[\psi\big(t^{n+1}\big)-\overline{\psi}^{n+1}\big]||^2+\epsilon k\nu||e_\omega^{n+1}||^2 \nonumber \\
		&\leq C(\epsilon, \nu, R_0, R_1) k \big[||D \omega \big(t^{n+1}\big)||^2+||\nabla^{-1}\eVO^{n}||_{G}^2\big] +\epsilon k\nu||e_\omega^{n+1}||^2,
	\end{align}
    where $\epsilon>0$ is a constant to be determined.
	Likewise,
	\begin{align}
		\label{non-Tra-2}
		&I_2\nonumber  \\
		&\leq C(\epsilon, \nu) k\left\{||D \omega\big(t^{n+1}\big)|| \cdot||\overline{\omega}^{n+1} || \cdot ||e_\omega^{n+1}|| +||\nabla \big[\nabla^{\perp}_H\overline{\psi}^{n+1}\cdot \nabla \Delta^{-1}\eo^{n+1}\big]|| \cdot ||\nabla^{-1}\eVO^{n}||_{G} \right\} \nonumber \\
		&\leq C(\epsilon, |\Omega|) k\left\{||D \omega\big(t^{n+1}\big)|| \cdot||\overline{\omega}^{n+1} ||   +||\nabla\overline{\omega}^{n+1}|| \cdot ||\nabla^{-1}\eVO^{n}||_{G} \right\} ||e_\omega^{n+1}|| \nonumber\\ 
		&\leq C(\epsilon, R_0, R_1, \nu, |\Omega|)k\big[||D \omega \big(t^{n+1}\big)||^2+||\nabla^{-1}\eVO^{n}||_{G}^2\big] +\epsilon k \nu||e_\omega^{n+1}||^2.
	\end{align}
	By the invariance of the global attractor under $\mathbb{S}_{q, k}$ and the consistency estimate \eqref{asym-con-eq}, one obtains
	\begin{align}
		\label{non-Tra-3}
		I_3 &\leq C(\epsilon, R_0, R_1, \nu, |\Omega|)k \big|q^{n+1}-1\big|^2 +\epsilon k\nu||e_\omega^{n+1}||^2 \nonumber \\
		&\leq C(\epsilon, R_0, R_1, \nu, |\Omega|)k^2+\epsilon k\nu||e_\omega^{n+1}||^2.
	\end{align}  
	For $I_4$,  by \eqref{diff-con1}  one has
	\begin{align}
		I_4 &\leq k\big(||D\psi_x\big(t^{n+1}\big)||+||[2e_\psi^n-e_{\psi}^{n-1}]_x||\big)\cdot ||\Delta^{-1}\eo^{n+1}|| \nonumber \\
		&\leq C(\epsilon, \nu)k\big[||D\omega\big(t^{n+1}\big)||^2+ ||\nabla^{-1}\eVO^{n}||_{G}^2\big] +\epsilon k\nu||e_\omega^{n+1}||^2. \label{non-Tra-32}
	\end{align}
	Finally one estimates $I_5$ as follows
	\begin{align}
		\label{non-Tra-4}
		I_5 &\leq k||\Delta^{-1} R_\omega^{n+1}|| \cdot ||e_\omega^{n+1}|| \nonumber \\
		&\leq C(\epsilon, \nu) k||\Delta^{-1} R_\omega^{n+1}||^2+\epsilon k\nu||e_\omega^{n+1}||^2.
	\end{align} 
	Taking $\epsilon=\frac{1}{6}$, and by virtue of \eqref{non-Tra-1}--\eqref{non-Tra-4}, one has
	\begin{align}
		\label{Err-SV-es-2}
		&||\nabla^{-1}\eVO^{n+1}||_{G}^2-||\nabla^{-1}\eVO^{n}||_{G}^2+\frac{\nu k}{6} ||\eVO^{n+1}||^2 \nonumber \\
		&\leq C(\epsilon, R_0, R_1, \nu, |\Omega|)k\big[||D \omega \big(t^{n+1}\big)||^2+||\nabla^{-1}\eVO^{n}||_{G}^2+||\Delta^{-1} R_\omega^{n+1}||^2+k\big].
	\end{align}  
	The consistency estimate \eqref{asym-con-eq} implies
	\begin{align*}
		||\nabla^{-1}\eVO^{1}||_{G} \leq C||\nabla^{-1}\big[\omega(t^1)-\omega^0\big]||+C||\nabla^{-1}\big[\omega^0-\omega^1\big]|| \leq k \sup_{t\in [0, T]} ||\nabla^{-1}\omega_t||+Ck.
	\end{align*}
    {  We remark that the free choice of $\omega^1$ contributes to the sub-optimal error estimate here. }
	Taking summation of \eqref{Err-SV-es-2} from $n=1$ to $m$, one arrives at
	\begin{align*}
		&||\nabla^{-1}\eVO^{m+1}||_{G}^2+\frac{\nu k}{6} \sum_{n=1}^m||\eVO^{n+1}||^2 \leq ||\nabla^{-1}\eVO^{1}||_{G}^2\nonumber \\
		&+C(R_0, R_1, \nu, |\Omega|)k\sum_{n=1}^m\Big[||D \omega \big(t^{n+1}\big)||^2+||\nabla^{-1}\eVO^{n}||_{G}^2+||\Delta^{-1} R_\omega^{n+1}||^2+k\Big], \nonumber \\
		&\leq Ck^2+C( R_0, R_1, \nu, |\Omega|, T)\left\{k\sum_{n=1}^m||\nabla^{-1}\eVO^{n}||_{G}^2+k^2\int_0^T||\omega_t||^2+||\Delta^{-1}\omega_{tt}||^2\, dt+k\right\}.
	\end{align*} 
	In light of the regularity of $\omega$ in \eqref{Reg-Tra}, an application of Gronwall's inequality  then gives
	\begin{align}\label{Tra-conv-Ne}
		||\nabla^{-1} [\omega\big(t^{n}\big)-\omega^{n}]||^2 \leq C(T, R_0, R_1, \nu, |\Omega|) k, \quad  t^n \leq T.
	\end{align}
	Through interpolation between $H^1_{p0}(\Omega)$ and $H^{-1}_{p}$, one concludes for any $\alpha \in [\frac{1}{4}, 1)$
	\begin{equation}\label{Tra-conv-alpha}
		\begin{split}
		&||\omega\big(t^{n}\big)-\omega^{n}||_{H^\alpha} \leq C ||\nabla^{-1} [\omega\big(t^{n}\big)-\omega^{n}]||^\frac{1-\alpha}{2} ||\nabla [\omega\big(t^{n}\big)-\omega^{n}]||^\frac{1+\alpha}{2} \\
		&\leq C(T, R_0, R_1, \nu, |\Omega|) k^{\frac{1-\alpha}{4}}.
		\end{split}
	\end{equation}
	This completes the proof of Lemma \ref{Traj-convL}.
\end{proof}

\bibliographystyle{amsplain}
\bibliography{multiphase-2025.bib}

\end{document}